\documentclass[10pt]{amsart}
\usepackage{amssymb,amsmath,latexsym,amsfonts}
\usepackage{tikz}
\usepackage{wasysym}
\usetikzlibrary{calc}

\def\ignore #1 {}

\newtheorem{thm}{Theorem}
\newtheorem{lem}[thm]{Lemma}
\newtheorem{defn}[thm]{Definition}

\newtheorem{cor}[thm]{Corollary}
\newtheorem{conjecture}{Conjecture}
\newtheorem{question}{Question}
\theoremstyle{definition}
\newtheorem*{remark}{Remark}
\newtheorem*{remarks}{Remarks}
\newtheorem{example}{Example}
\newtheorem{algorithm}[thm]{Algorithm}
\newtheorem*{ut}{Useful trick}

\def\hpic #1 #2 {\mbox{$\begin{array}[c]{l} \epsfig{file=#1,height=#2} \end{arr\
ay}$}}
\def\vpic #1 #2 {\mbox{$\begin{array}[c]{l} \epsfig{file=#1,width=#2} \end{arra\
y}$}}

\def\R{\mbox{{\bf R}}}
\def\C{\mbox{{\bf C}}}

\def\a{\mbox{${\alpha}$}}

\newcommand{\bdry}{\partial}

\DeclareMathOperator{\im}{Im}

\DeclareMathOperator{\Aff}{Aff}
\DeclareMathOperator{\Area}{Area}
\DeclareMathOperator{\Vxs}{Vertices}
\DeclareMathOperator{\Triangles}{Triangles}
\DeclareMathOperator{\codim}{codim}

\newcommand{\D}{\Delta}
\renewcommand{\a}{\alpha}

\newcommand{\s}{\sigma}

\newcommand{\T}{{\mathcal{T}}}

\def \sb #1 {\overline{\s}_{#1}}
\def\l{L}

\begin{document}
\title{Spaces of polygonal triangulations and Monsky polynomials\footnote{This article was published in \emph{Discrete and Computational Geometry}, vol.~51 no.~1 (2014), pp.~132--160.}}
\author{Aaron Abrams and James Pommersheim}

\address{
Mathematics Department, Washington and Lee University, Lexington VA 24450 \newline
\indent
Department of Mathematics, Reed College, 3203 SE Woodstock Blvd, Portland OR 97202}


\begin{abstract} Given a combinatorial triangulation of an $n$-gon, 
we study (a) the space of all possible drawings in the plane such the edges are straight line segments and 
the boundary has
a fixed shape, and (b) the algebraic variety 
of possibilities
for the areas of the triangles in such drawings.  We define a generalized notion of 
triangulation, and we show that the areas of the triangles in a generalized
triangulation $\T$ of a square must satisfy a single irreducible homogeneous 
polynomial relation $p(\T)$ depending only on the combinatorics of $\T$.  
The invariant $p(\T)$ is called the \emph{Monsky polynomial}; it 
captures algebraic, geometric, and combinatorial information about $\T$.  We 
give an algorithm that computes a lower bound on
the degree of $p(\T)$, and we present several examples in which the
algorithm is used to compute the degree.
\end{abstract}
%
%
\dedicatory{The first author dedicates this work to Sherman K.~Stein, from whom he first 
learned about Monsky's Theorem.}

\maketitle

\tableofcontents

\section{Introduction}
\subsection{Background}
This paper grew out of a desire to better understand a theorem of Paul Monsky
about plane geometry.  Monsky's theorem from 1970 (\cite{monsky}, 
see also \cite{steinbook}) states that it
is impossible to cut a square into an odd number of triangles, all of the same area.
This result grew out of a problem that was first posed in the 
American Mathematical Monthly by Fred Richman and John Thomas
\cite{richman} and later advanced by Thomas \cite{thomas}.
Monsky's proof blends combinatorics, topology, and number theory in an elementary 
but delicate way.

Monsky's result has generally been viewed as a statement about ``equidissections,''
or tilings by equal area tiles.  As such it has been generalized in various directions, 
e.g.~to polygons other than a square \cite{stein1,stein2,monsky2} and to higher dimensional 
cubes \cite{mead}.  These results all have a combinatorial feel to them, and the
rigidity of the notion of an equidissection is crucial.

With the hope of understanding Monsky's theorem in a different way, 
Joe Buhler suggested studying, for a fixed combinatorial triangulation of
a square, the space of all realizable collections of areas of triangles. 

\begin{question}\label{buhler}
Given the combinatorics of a triangulation $\T$ of the unit square, which collections
of areas can be realized by the triangles?
\end{question}

A complete answer to this question would, of course, include Monsky's theorem.
One observes that the areas of the triangles add up to the area of the square, and by a
dimension count, one expects the areas of the triangles to satisfy exactly one
additional polynomial equation.  In his senior thesis,
Buhler's student Adam Robins calculated this polynomial $p(\T)$for several
particular triangulations \cite{transcendence}.

\begin{example}\label{example1}
If $\T$ is the triangulation whose combinatorial type is drawn in Figure \ref{fig:X}, 
then no matter where the central vertex is placed, the four areas must satisfy 
$p(\T)=0$, where $p(\T):=A-B+C-D$.

\begin{figure}[ht]
\begin{center}
\begin{tikzpicture}
\draw (0,0) -- (0,3) -- (3,3) -- (3,0) -- cycle;
\draw (0,0) -- (1.8,2.1) -- (3,3);
\draw (0,3) -- (1.8,2.1) -- (3,0);
\draw (1.5,0.8) node{$A$};
\draw (0.7,1.9) node{$B$};
\draw (2.5,1.9) node{$D$};
\draw (1.7,2.6) node{$C$};
\end{tikzpicture}
\caption{In any drawing of this 2-complex, as long as the boundary is
a parallelogram the areas will satisfy $A-B+C-D=0$, regardless of the placement
of the central vertex.}
\label{fig:X}
\end{center}
\end{figure}
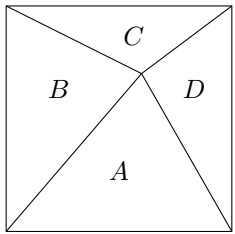
\end{example}

A few years after \cite{transcendence}, in an
REU supervised by Serge Tabachnikov and Misha Guysinsky, 
students Josh Kantor and Max Maydanskiy obtained further results about
this polynomial, including the calculation of several more examples 
and proofs of some general bounds on its degree \cite{gonewild}.
These authors also give a family of examples showing that the degree
can be exponential in the size of the triangulation.

Neither \cite{transcendence} nor \cite{gonewild} gives a rigorous proof
of existence of the polynomial they study.  Monsky's original theorem from
\cite{monsky} does prove that the areas satisfy a polynomial relation, although
his approach was quite different from ours and it was initially unclear to us whether 
his polynomial was the same as ours.  In recent communication with Monsky,
however, it has become clear that his and ours are versions of the same thing.  Thus
in Section \ref{sec:Monsky} we name the polynomial the Monsky polynomial.

Our approach to Monsky's polynomial is the following.  We study 
(a) the space of all possible drawings of a given abstract triangulation $\T$ of an 
$n$-gon in which the boundary has a fixed shape, and (b) the 
variety $V(\T)$ of possibilities for the areas of the triangles in such drawings.  
We establish a framework allows the possibility of degenerating triangles by
requiring that certain subsets of the vertices be collinear.
For convenience we work over the complex numbers; in this framework,
we show that $V(\T)$ is an irreducible algebraic variety, and we compute its dimension. 
In the case $n=4$, $V(\T)$ is a projective hypersurface, and so there is a single irreducible polynomial
relation $p(\T)$ that is satisfied by the areas of the triangles.  This is the Monsky polynomial.

\subsection{Results}
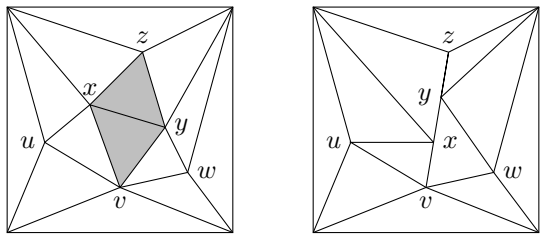
\begin{figure}[ht]
\begin{center}
\begin{tikzpicture}[-]

\matrix (m) [row sep = 8em, column sep = 3em]
{
\coordinate (P) at (0,0);
\coordinate (Q) at (3,0);
\coordinate (R) at (3,3);
\coordinate (S) at (0,3);

\draw (P) -- (Q) -- (R) -- (S) -- cycle;

\coordinate (u) at (0.5,1.2);
\coordinate (v) at (1.5,0.6);
\coordinate (w) at (2.4,0.8);
\coordinate (x) at (1.1,1.7);
\coordinate (y) at (2.1,1.4);
\coordinate (z) at (1.8,2.4);

\fill[gray!50!white] (x)--(y)--(z)--cycle;
\fill[gray!50!white] (x)--(y)--(v)--cycle;

\draw (u) node[anchor=east]{$u$}--
	(v) node[anchor=north]{$v$} --
	(w) node[anchor=west]{$w$}--
	(y) node[anchor=west]{$y$}--
	(z) node[anchor=south]{$z$}--
	(x) node[anchor=south]{$x$}--
	cycle;
\draw (v)--(x)--(y)--cycle;
\draw (P) -- (u);
\draw (P) -- (v);
\draw (S) -- (u);
\draw (S) -- (x);
\draw (S) -- (z);
\draw (Q) -- (v);
\draw (Q) -- (w);
\draw (R) -- (w);
\draw (R) -- (y);
\draw (R) -- (z);


&

\draw (P) -- (Q) -- (R) -- (S) -- cycle;

\coordinate (u) at (0.5,1.2);
\coordinate (v) at (1.5,0.6);
\coordinate (w) at (2.4,0.8);
\coordinate (x) at (1.6,1.2);
\coordinate (y) at (1.7,1.8);
\coordinate (z) at (1.8,2.4);

\draw (u) node[anchor=east]{$u$}--
	(v) node[anchor=north]{$v$} --
	(w) node[anchor=west]{$w$}--
	(y) node[anchor=east]{$y$}--
	(z) node[anchor=south]{$z$}--
	(x) node[anchor=west]{$x$}--
	cycle;
\draw (v)--(x);
\draw (P) -- (u);
\draw (P) -- (v);
\draw (S) -- (u);
\draw (S) -- (x);
\draw (S) -- (z);
\draw (Q) -- (v);
\draw (Q) -- (w);
\draw (R) -- (w);
\draw (R) -- (y);
\draw (R) -- (z);
\\};

\end{tikzpicture}
\caption{On the left is a depiction of an abstract triangulation $\T=(T,C)$.  The shading
indicates the collinearity condition $C=\{v,x,y,z\}$, which is the set of vertices of the shaded 
triangles.  On the right is (the image of) a typical drawing of the triangulation $(T,C)$.  The 
shaded triangles have degenerated to segments because the vertices of $C$ are collinear. }
\label{fig:drawing}
\end{center}
\end{figure}

We now give an informal summary of the results of this paper.
We begin with a combinatorial triangulation $T$ of an $n$-gon.
We define a \emph{triangulation} $\T$ to be a pair $(T,C)$, where
$C$ is the set of vertices of a contiguous subset of triangles of $T$. 
(See Figure \ref{fig:drawing} for an example, and Section \ref{sec:triangulations} 
for the precise definitions.) Intuitively we think of $C$ as
a ``collinearity condition,'' i.e., a set of vertices which will
be required to be collinear.  A drawing of $\T$ is a placement of
the vertices of $T$ in the plane with the condition that the vertices of $C$ must be
collinear.   We then think of the edges of the triangulation being drawn as straight line segments.
We consider the space $X(\T)$ of drawings of $\T$, and in Theorem
\ref{irr} we show
that this space is an irreducible algebraic variety.  

We next consider the space of drawings of $\T$ in which the
boundary is constrained to be a fixed shape.  More precisely, if 
$\pentagon$ is an $n$-gon in the plane, we define $X_{\pentagon}(\T)$
to be the subset of $X(\T)$ consisting of those drawings in which
the boundary of $T$ is mapped to an affine image of $\pentagon$.
In Theorem \ref{irr2} we
show that $X_{\pentagon}(\T)$ is also an irreducible algebraic variety.

Next, we introduce the area map $f_{\pentagon}(\T)$ corresponding to a triangulation
$\T$ and a boundary shape $\pentagon$.  To do this we let $Y(\T)$ be the
projective space with one coordinate
for each non-degenerate triangle of $\T$.  Then the rational map 
$$f_{\pentagon}(\T):X_{\pentagon}(\T) \dashrightarrow Y(\T)$$
records the areas of the triangles at each point of $X_{\pentagon}(\T)$.  The 
closure of the image of $f_{\pentagon}(\T)$ is the variety $V_{\pentagon}(\T)$ of possible
areas.  In Theorem \ref{thm:dimension} we show that under suitable
hypotheses, $V_{\pentagon}(\T)$ has codimension $n-3$ in $Y(\T)$.

When $n=4$, so that $T$ is a quadrilateral, Theorem \ref{thm:dimension} 
says that $V_{\pentagon}(\T)$ has codimension 1 in $Y(\T)$.  Of special interest
is the case that the boundary shape $\pentagon$ is a square.  In this case
we define $p(\T)$ to be the irreducible polynomial that vanishes on 
$V_{\pentagon}(\T)$.  The polynomial $p(\T)$ is the \emph{Monsky polynomial}.

The framework we establish allows us to successively degenerate triangles
to produce new triangulations.
In Section \ref{sec:degree} we use this feature to give a purely combinatorial algorithm that 
recursively computes a lower bound for $\deg(\T):=\deg(p(\T))$.  
The algorithm exploits the interplay between the combinatorics 
of $\T$, the geometry of $V(\T)$, and the algebra of $p(\T)$.  

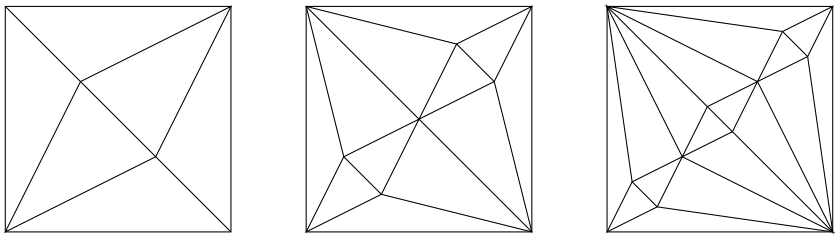
\begin{figure}[ht]
\begin{center}
\begin{tikzpicture}
\draw (0,0) -- (0,3) -- (3,3) -- (3,0) -- cycle;
\draw (0,0) -- (1,2) -- (3,3) -- (2,1) -- cycle;
\draw (3,0) -- (0,3);

\draw (4,0) -- (4,3) -- (7,3) -- (7,0) -- cycle;
\draw (4,0) -- (4.5,1) -- (5.5,1.5) -- (5,0.5) -- cycle;
\draw (5.5,1.5) -- (6,2.5) -- (7,3) -- (6.5,2) -- cycle;
\draw (4,3) -- (4.5,1) -- (5,0.5) -- (7,0) -- (6.5,2) -- (6,2.5) -- cycle;
\draw (4,3) -- (7,0);

\draw (8,0) -- (8,3) -- (11,3) -- (11,0) -- cycle;
\draw (8,0) -- (25/3,2/3) -- (9,1) -- (26/3,1/3) -- cycle;
\draw (9,1) -- (28/3,5/3) -- (10,2) -- (29/3,4/3) -- cycle;
\draw (10,2) -- (31/3,8/3) -- (11,3) -- (32/3,7/3) -- cycle;
\draw (8,3) -- (25/3,2/3) -- (26/3,1/3) -- (11,0) -- (9,1) -- cycle;
\draw (8,3) -- (10,2) -- (11,0) -- (32/3,7/3) -- (31/3,8/3) -- cycle;
\draw (8,3) -- (11,0);

\end{tikzpicture}
\caption{The degree of $p(\T)$ for the $k$th triangulation in this sequence
is exponential in $k$.  Conjecturally, the degree is exactly $(3k-2)2^{k-1}$ for
$k>1$.  See Section \ref{sec:examples}.}
\label{fig:exponential}
\end{center}
\end{figure}

As consequences of the algorithm,
we characterize those triangulations $\T$ for which the degree of $p(\T)$
is $1$, and we prove that if there is no collinearity condition ($C=\emptyset$) and $T$ cannot be
obtained by subdividing a simpler triangulation, then  
the degree of $p(\T)$ is bounded below by the number of interior vertices of $T$.
Then, in Section \ref{sec:examples} we give an infinite family of examples for which 
our algorithm computes the correct degree of $\T$, as well as an infinite family
for which the degree grows exponentially in the number of vertices of $T$.
See Figure \ref{fig:exponential}.

\subsection{Further remarks}
It is worth mentioning that for any particular $\T$, the polynomial $p(\T)$ is theoretically 
computable:  to do so is an instance of the well-known \emph{implicitization
problem} in algebraic geometry, whose solution involves techniques of
Gr\"obner bases or resultants (see \cite{CLO, MC}).  Unfortunately, using these
general techniques to determine $p(\T)$ is intractable in all but 
the smallest cases.  Even when feasible, these computations give little insight into 
the general structure of the polynomial $p(\T)$.  In contrast, using code written
by Phillip Andreae \cite{phillip} for his senior thesis, we were able
to run Algorithm \ref{alg} on hundreds of triangulations with a dozen or more
vertices.

Our work has a similar flavor to recent work on Robbins' 
conjectures (\cite{robbins1,robbins2,pak,paksurvey,connelly,varfolomeev}),
in which a polynomial relation among squared lengths of the sides of a cyclic polygon is 
studied (a special case is Heron's formula).
There also appear to be similarities with the well-known Bellows Conjecture 
(\cite{sabitov,connelly2,csw}).

\subsection{Acknowledgements}
Thanks to Dave Perkinson, Joe Buhler, Serge Tabachnikov, Phillip Andreae, Joe Fu,
and Ray Mayer for numerous stimulating conversations on this subject.  
Special thanks to Dave Perkinson for introducing us to this problem
and for writing an early CoCoA script to compute several examples.
Thanks to Phillip Andreae for producing a senior thesis \cite{phillip} that included
java code to implement our algorithm before the algorithm made any sense.  
Finally, thanks to the Research in Pairs program at MFO (Oberwolfach) 
for allowing us to spend two weeks in 2011 converting ham into theorems.

\section{Triangulations of polygons} \label{sec:triangulations}

Consider a finite 2-dimensional simplicial complex $T$ whose 
underlying space $|T|$ is homeomorphic to a disk.  Let $\Vxs(T)$ denote the set of vertices of $T$.
We orient the complex and label the boundary vertices 
$P_1,\ldots,P_n$ in cyclic order according to the orientation.  The $P_i$
are called the \emph{corners} of $T$.  The set of corners
is denoted by $\bdry T$.  Vertices of $T$ that are not corners are
called \emph{interior vertices}.

It is an easy exercise to show that the number of triangles in $T$
is $2k+n-2$, where $k$ is the number of interior vertices.
If $S$ is a set of triangles of $T$ then we
define a graph $S^*$ with one vertex for each triangle in $S$ 
and an edge connecting pairs of adjacent triangles of $S$.  
(Note that $S^*$ is just the subgraph of the dual graph to the triangulation $T$
spanned by the vertices that correspond to $S$.) 
We say that $S$ is \emph{contiguous} if $S^*$ is a connected graph.
(Equivalently, $S$ is contiguous if the interior of the union of the (closed)
triangles in $S$ is connected.)
We denote by $\Vxs(S)$ the set
of vertices $v\in\Vxs(T)$ such that $v$ is a vertex of a triangle of $S$.

\begin{defn}[Condition]\label{def:condition}
A \emph{condition} $C$ is a collection of vertices of $T$ such that
\begin{itemize}
\item 
$C=\Vxs(S)$ for some contiguous set of triangles $S$, and
\item
$|C\cap\bdry T| \le 2$.
\end{itemize}
\end{defn}

\begin{figure}[ht]
\begin{center}
\begin{tikzpicture}[-]

\matrix (m) [row sep = 8em, column sep = 3em]
{
\coordinate (P) at (0,0);
\coordinate (Q) at (3,0);
\coordinate (R) at (3,3);
\coordinate (S) at (0,3);

\coordinate (t) at (0.6,2.2);
\coordinate (u) at (0.5,1.2);
\coordinate (v) at (1.5,0.6);
\coordinate (w) at (2.4,0.8);
\coordinate (x) at (1.1,1.7);
\coordinate (y) at (2.1,1.4);
\coordinate (z) at (1.8,2.4);

\coordinate (A) at (1.6,1.2);
\coordinate (B) at (1.7,1.8);
\coordinate (C) at (1,1.2);
\coordinate (D) at (0.6,0.6);
\coordinate (E) at (0.25,1.2);

\fill[gray!50!white] (P)--(S)--(u)--cycle;
\fill[gray!50!white] (P)--(v)--(u)--cycle;
\fill[gray!50!white] (x)--(v)--(u)--cycle;
\fill[gray!50!white] (x)--(v)--(y)--cycle;
\fill[gray!50!white] (x)--(z)--(y)--cycle;

\draw (P) -- (Q) -- (R) -- (S) -- cycle;

\draw (u) --(v) --(w) --(y) --(z)--(x) -- (u)--(t)--(z);
\draw (t) -- (x);
\draw (v)--(x)--(y)--cycle;
\foreach \p in {u,v} \draw (P) -- (\p);
\foreach \p in {t,u,z} \draw (S) -- (\p);
\foreach \p in {v,w} \draw (Q) -- (\p);
\foreach \p in {w,y,z} \draw (R) -- (\p);

\draw[fill=black] (A) circle (1pt);
\draw[fill=black] (B) circle (1pt);
\draw[fill=black] (C) circle (1pt);
\draw[fill=black] (D) circle (1pt);
\draw[fill=black] (E) circle (1pt);

\draw[dotted] (B) -- (A) -- (C) -- (D) -- (E);

&

\fill[gray!50!white] (P)--(S)--(u)--cycle;
\fill[gray!50!white] (P)--(v)--(u)--cycle;
\fill[gray!50!white] (y)--(v)--(w)--cycle;

\draw (P) -- (Q) -- (R) -- (S) -- cycle;

\draw (u) --(v) --(w) --(y) --(z)--(x) -- (u)--(t)--(z);
\draw (t) -- (x);
\draw (v)--(x)--(y)--cycle;
\foreach \p in {u,v} \draw (P) -- (\p);
\foreach \p in {t,u,z} \draw (S) -- (\p);
\foreach \p in {v,w} \draw (Q) -- (\p);
\foreach \p in {w,y,z} \draw (R) -- (\p);

\coordinate (F) at (2.05,1);
\coordinate (G) at (2.3,0.5);

\draw[fill=black] (D) circle (1pt);
\draw[fill=black] (E) circle (1pt);
\draw[fill=black] (F) circle (1pt);

\draw[dotted] (E)--(D);

&

\fill[gray!50!white] (P)--(S)--(u)--cycle;
\fill[gray!50!white] (P)--(v)--(u)--cycle;
\fill[gray!50!white] (y)--(v)--(x)--cycle;

\draw (P) -- (Q) -- (R) -- (S) -- cycle;

\draw (u) --(v) --(w) --(y) --(z)--(x) -- (u)--(t)--(z);
\draw (t) -- (x);
\draw (v)--(x)--(y)--cycle;
\foreach \p in {u,v} \draw (P) -- (\p);
\foreach \p in {t,u,z} \draw (S) -- (\p);
\foreach \p in {v,w} \draw (Q) -- (\p);
\foreach \p in {w,y,z} \draw (R) -- (\p);

\coordinate (F) at (2.05,1);
\coordinate (G) at (2.3,0.5);

\draw[fill=black] (A) circle (1pt);
\draw[fill=black] (D) circle (1pt);
\draw[fill=black] (E) circle (1pt);

\draw[dotted] (E)--(D);

\draw (C) node {$\Delta$};
\\};

\end{tikzpicture}
\caption{In each of these figures, let $S$ denote the set of shaded triangles.  The
dual graph $S^*$ is also shown with dotted lines.
On the left, $S$ is contiguous and $\Vxs(S)$ is a condition.  In the middle, $S$ is
not contiguous and $\Vxs(S)$ is not a condition.  On the right, $S$ is not contiguous
but $\Vxs(S)$ is a condition because $S'=S\cup \Delta$ is contiguous and $\Vxs(S')=\Vxs(S)$.
Note $S'=\Triangles(\Vxs(S))$.
}
\label{fig:conditions}
\end{center}
\end{figure}
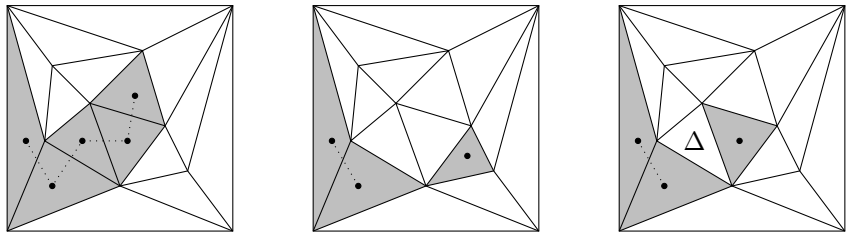

See Figure \ref{fig:conditions} for examples.
We think of $C$ as a collection of vertices that will be required to be collinear; 
$C$ arises from ``killing" a connected region of triangles.

For a given condition $C$, note that there may be more than one $S$ such that
$C=\Vxs(S)$, but there is a uniquely determined maximal such $S$ consisting of
all those triangles whose vertices are all in $C$.  We denote this set of triangles
by $\Triangles(C)$.  Again see Figure \ref{fig:conditions}.

If $C$ is a condition then we abbreviate $(\Triangles(C))^*$ by $C^*$.

A \emph{triangulation} $\T$ is a pair $(T,C)$ where $T$ is a simplicial complex
as described above, and $C$ is a (possibly empty) condition. 

\section{Drawings of triangulations}\label{sec:X}

We are interested in drawing triangulations $(T,C)$ in the affine plane.  We will
study all drawings of $(T,C)$, first without restricting the boundary and then
later with the constraint that the boundary forms a fixed shape, defined up to
affine equivalence.

We will work over the field $\C$ of complex numbers,
but many of the arguments are more easily visualized over the reals.

A fundamental object associated to a triangulation $\T=(T,C)$ is the domain
$$X=X(\T)$$
defined to be the set of those maps 
$\rho:\Vxs(T)\to \C^2$ 
such that there is an affine line $\ell\subset \C^2$ containing $\rho(v)$ for every
vertex $v$ of the condition $C$.  Note that $X(\T)$ is an affine subvariety of the space
of all maps $\Vxs(T)\to\C^2$.  Each point $\rho\in X(\T)$ is called a
\emph{drawing} of $\T$.

\begin{thm}\label{irr}
If $\T$ is a triangulation, then 
$X(\T)$ is an irreducible algebraic variety over $\C$.
\end{thm}

\begin{proof}

Denote the interior vertices of $T$ by $v_1,\ldots,v_k$.
Consider any total order $<$ on $\Vxs(T)$ satisfying 
$P_i<v_j$ for all $i$ and $j$.
For $v\in \Vxs(T)$ let $\a(v)=1$ if $v\in C$ and there are
at least two vertices $x_1,x_2\in C$ with $x_i<v$.  Otherwise
let $\a(v)=2$.  Then, if we draw $\T$ by placing the vertices
in order, then $\a(v)$ is the number of degrees of freedom
we have in placing $v$.

We parameterize $X$ by defining a rational map
$$g:\prod\limits_{v\in\Vxs(T)} \C^{\a(v)} \longrightarrow X.$$
Let $x=(x_v)_{v\in\Vxs(T)}$ be coordinates on the domain 
$\prod\limits_{v\in\Vxs(T)} \C^{\a_v} $.
The coordinate functions $\rho(v)(x)$ of the point $g(x)$ are constructed 
inductively, as follows.
For each $v\in\Vxs(T)$, we assume that the 
coordinate functions $\rho(w)(x)$ for vertices $w$ with $w<v$ have been defined.
Then, if $\a(v)=2$, we set $\rho(v)(x)=x_v$.  If $\a(v)=1$, then set 
$\rho(v)(x)=x_v \rho(y)(x) + (1-x_v) \rho(z)(x)$
where $y$ and $z$ are the first two points of $C$ and $y<z$.

Clearly the image of $g$ is contained in $X$.
It then follows that $X$ is irreducible, provided that the 
image of $g$ is a Zariski dense subset of $X$.  To show this, take any drawing 
$\rho\in X$.  Note that $\rho$ is in the image of $g$ iff either $\rho(y)\ne\rho(z)$
or all the points $\rho(v)$ coincide for $v\in C$.  Hence if $\rho$ is not in the image
of $g$, then $\rho(y)=\rho(z)$ but the line $\ell$ is uniquely defined.  By perturbing
$\rho(z)$ along $\ell$ we obtain a nearby point $\rho'$ of $X$ that is in the image
of $g$.  Thus the image of $g$ is dense in $X$, and $X$ is irreducible.
\end{proof}

The next result is similar, but we impose a boundary condition.
Given $n\ge 3$, an \emph{$n$-gon shape} $\pentagon$ is an ordered $n$-tuple 
$(z_1,\ldots,z_n)$ of points in $\C^2$ such that no three are collinear.
Suppose $\T=(T,C)$ is given with $|\bdry T|=n$, and let $\pentagon$
be an $n$-gon shape.  We let $X_{\pentagon}(\T)$ be 
the set of drawings of $\T$ whose boundary is an affine image of $\pentagon$;
that is,
$$X_{\pentagon}(\T)=\{\rho\in X(\T) \mid \exists \phi\in\widehat\Aff \mbox{ with }
\rho(P_i)=\phi(z_i) \ \forall i\},$$
where $\widehat\Aff$ is the set of (not necessarily invertible) affine transformations
of $\C^2$.

\begin{thm}\label{irr2}
If $\T$ is a triangulation with $n$ boundary vertices and $\pentagon$ is
an $n$-gon shape, then $X_{\pentagon}(\T)$ is an irreducible algebraic variety over $\C$.
\end{thm}

\begin{proof}
We proceed basically as in the previous theorem, with some 
modifications.  We order the vertices in any way as long as $P_i<v_j$ for all
$i$ and $j$ and the vertices of $C\cap \bdry T$ come first.
Set $\a(v)=2$ for the first three vertices in the order, set $\a(v)=0$ for
any remaining boundary vertices, and define $\a(v_i)$ as before.

We parameterize $X_{\pentagon}$ again by $\prod\limits_{v\in\Vxs(T)} \C^{\a(v)}$.
Drawing the vertices in order,
we see that the first three points determine
the rest of the boundary as well as the affine map $\phi$.
For the $v_i$, the parameterization works as in Theorem \ref{irr}.

Again we wish to show that the image of the parameterization is 
dense in $X_{\pentagon}(\T)$.  Points $\rho\in X_{\pentagon}(\T)$
in the image of the parameterization are characterized exactly
as before, and if a point $\rho$ of $X$ is not in the image then
by our choice of ordering on vertices, the second point of the
condition has $\a=2$, and so we can move it along $\ell$ as
in the previous proof.
\end{proof}

\begin{remarks}
The above theorems easily generalize to triangulations with several
collinearity conditions, as long as the conditions are disjoint.
The hypotheses on the condition $C$ are also unnecessary
for irreducibility; one could allow $C$ to be an arbitrary subset of $\Vxs(T)$.

With some additional arguments, the above theorems also generalize to so-called
\emph{dissections}, which is the context in which Monsky made his original
discoveries.  A dissection is a decomposition of the polygon into closed
triangles whose interiors are disjoint; the decomposition need not be
simplicial.  Every dissection can be viewed as a (simplicial) triangulation with
possibly many collinearity conditions.

Perhaps the most general setting of all is that of triangulations with arbitrary
collinearity conditions that are allowed to intersect.  One difficulty that
arises in this situation is that the conditions may not be logically 
independent, meaning that it is possible for one to follow from the others.  For
example, this could 
be the result of an incidence theorem such as Pappus's Theorem.  In such a case
$X$ may not have the expected dimension and it becomes more difficult
to parameterize.

Another difficulty with this generalization is that it can cause $X$ to
be reducible.  This happens, for instance, if $X$ contains one condition
of the form $\{w,x,y\}$ and another condition of the form $\{x,y,z\}$.
In this case $X$ will contain two components, one consisting of
drawings in which $w,x,y,z$ are all drawn collinearly and one 
consisting of drawings in which $x$ and $y$ are drawn copointally.
See Example \ref{example2}.

We believe that if every pair of collinearity conditions 
intersects in at most one vertex, then the space of drawings is irreducible.
\end{remarks}

\section{The variety of areas}\label{sec:areas}

If $\D$ is an oriented triangle in $\C^2$ with vertices $(x_1,y_1),(x_2,y_2),(x_3,y_3)$
then the \emph{area} of $\D$ is given by the formula
$$\Area(\D)=\frac 12
\left|
\begin{matrix}
1 & 1 & 1 \\
x_1 & x_2 & x_3 \\
y_1 & y_2 & y_3 
\end{matrix}
\right|.
$$
This is the usual (signed) area when $\D\subset\R^2$.

Let the triangulation $\T$ be given, let $n$ be the number of boundary points, 
and let $\pentagon$ be a fixed $n$-gon shape.
Let $Y=Y(\T)$ be the projective space with a coordinate 
$A_i$ for each triangle $\D_i$ whose vertices are not all in $C$.   We have rational maps
$f=f(\T):X(\T)\dashrightarrow Y(\T)$
and $f_{\pentagon}(\T):X_{\pentagon}(\T)\dashrightarrow Y(\T)$
called the \emph{area maps} given by 
$$f_*(\rho)=[\cdots:\Area(\rho(\D_i)):\cdots].$$

The main object of study in this paper is the closure of the image of $f_{\pentagon}$,
denoted by $V_{\pentagon}=V_{\pentagon}(\T)$.  This is an
irreducible (projective) subvariety of $Y$.

The goal of this section is to establish the dimension of $V_{\pentagon}$.
For a point $\rho$ of $X_{\pentagon}$ with $\rho(P_1),\rho(P_2),\rho(P_3)$ not
collinear, the orbit of $\rho$ under the action of the group $\Aff$ of invertible affine
transformations of $\C^2$ is contained in a single fiber of $f_{\pentagon}$.  Thus
generic fibers have dimension at least 6 (the dimension of $\Aff$).

\begin{figure}[ht]
\begin{center}
\begin{tikzpicture}[-]

\matrix (m) [row sep = 8em, column sep = 5em]
{

\coordinate (P) at (0,0);
\coordinate (Q) at (3,0);
\coordinate (R) at (3,3);
\coordinate (S) at (0,3);
\coordinate (1) at (.7,.8);
\coordinate (2) at (1.5,.7);
\coordinate (3) at (2.2,.7);
\coordinate (4) at (2.2,1.5);
\coordinate (5) at (2.2,2.3);
\coordinate (6) at (1.5,2.3);
\coordinate (7) at (.8,2.2);
\coordinate (8) at (.7,1.5);
\coordinate (9) at (1.2,1.2);
\coordinate (10) at (1.8,1.2);
\coordinate (11) at (1.8,1.8);
\coordinate (12) at (1.2,1.9);
\coordinate (13) at (1.45,1.55);

\foreach \x in {1,2,8}
	\draw[gray] (P)--(\x);
\foreach \x in {2,3,4}
	\draw[gray] (Q)--(\x);
\foreach \x in {4,5,6}
	\draw [gray](R)--(\x);
\foreach \x in {6,7,8}
	\draw[gray] (S)--(\x);
	
\fill[gray!50!white] (1)--(2)--(9)--cycle;
\fill[gray!50!white] (10)--(2)--(9)--cycle;
\fill[gray!50!white] (10)--(2)--(3)--cycle;
\fill[gray!50!white] (10)--(4)--(3)--cycle;
\fill[gray!50!white] (10)--(4)--(11)--cycle;
\fill[gray!50!white] (5)--(4)--(11)--cycle;
\fill[gray!50!white] (5)--(6)--(11)--cycle;
\fill[gray!50!white] (12)--(6)--(11)--cycle;
\fill[gray!50!white] (12)--(6)--(7)--cycle;
\fill[gray!50!white] (12)--(8)--(7)--cycle;
\fill[gray!50!white] (12)--(8)--(9)--cycle;
\fill[gray!50!white] (1)--(8)--(9)--cycle;

\draw (1)--(2)--(3)--(4)--(5)--(6)--(7)--(8)--cycle;
\draw (9)--(10)--(11)--(12)--cycle;
\foreach \x in {9,10,11,12}
	\draw[gray] (\x)--(13);
\draw (1)--(9)--(2)--(10)--(3);
\draw (10)--(4)--(11)--(5);
\draw (11)--(6)--(12)--(7);
\draw (12)--(8)--(9)--(1);
\draw (13) node[anchor=north]{$v$};

\draw (P)--(Q)--(R)--(S)--cycle;

&

\coordinate (1) at (1.4,.7);
\coordinate (2) at (1.4,.3);
\coordinate (3) at (1.4,.5);
\coordinate (4) at (1.4,1.3);
\coordinate (5) at (1.4,2.1);
\coordinate (6) at (1.4,2.5);
\coordinate (7) at (1.4,2.3);
\coordinate (8) at (1.4,1.5);
\coordinate (9) at (1.4,.9);
\coordinate (10) at (1.4,1.1);
\coordinate (11) at (1.4,1.7);
\coordinate (12) at (1.4,1.9);
\coordinate (13) at (2.3,1.55);

\foreach \x in {1,2,8}
	\draw[gray] (P)--(\x);
\foreach \x in {2,3,4}
	\draw[gray] (Q)--(\x);
\foreach \x in {4,5,6}
	\draw [gray](R)--(\x);
\foreach \x in {6,7,8}
	\draw[gray] (S)--(\x);

\draw (1)--(2)--(3)--(4)--(5)--(6)--(7)--(8)--cycle;
\draw (9)--(10)--(11)--(12)--cycle;
\foreach \x in {9,10,11,12}
	\draw[gray] (\x)--(13);
\draw (1)--(9)--(2)--(10)--(3);
\draw (10)--(4)--(11)--(5);
\draw (11)--(6)--(12)--(7);
\draw (12)--(8)--(9)--(1);
\fill[black] (13) node[anchor=west]{$v$} circle (1pt);

\draw (P)--(Q)--(R)--(S)--cycle;

\draw[->] (2.3,1.65)--(2.3,1.8);
\draw[->] (2.3,1.45)--(2.3,1.3);

\\};

\end{tikzpicture}
\caption{If $C$ consists of the vertices of the shaded triangles, then $C$ is a condition that
isn't non-separating.  In any drawing of $T$, there is a line containing all the vertices of $C$,
and the vertex $v$ can move parallel to this line without changing the area of any
triangle.}
\label{fig:sep}
\end{center}
\end{figure}
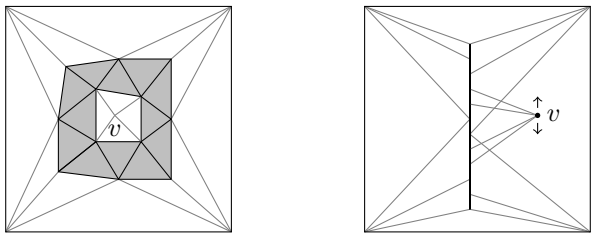

Exactly two phenomena can cause fibers to be larger than this.
First, consider the subgraph $G$ of the (1-skeleton of) $T$ spanned by 
those vertices that are not in $C$.  If this graph has a connected component
that contains none of the $P_i$, then this component can be 
``sheared" with respect to $\ell$, independently of the rest of $G$, 
without changing the area of any triangle.
Thus there is at least one additional 1-parameter family of drawings 
inside the fiber for each such component.  See Figure \ref{fig:sep}.

\begin{defn} \label{def:non-sep}[Non-separating]
Let $\T=(T,C)$.  Let $G$ be the graph consisting of those edges of $T$ whose
both endpoints are not in $C$.  The condition $C$ is \emph{non-separating} if
every connected component of $G$ contains at least one corner of $T$.
\end{defn}

(Recall that a corner is a vertex of $\bdry T$.)  Note that it is possible to have
$C$ non-separating and $G$ disconnected by having
$C$ contain two non-adjacent corners.

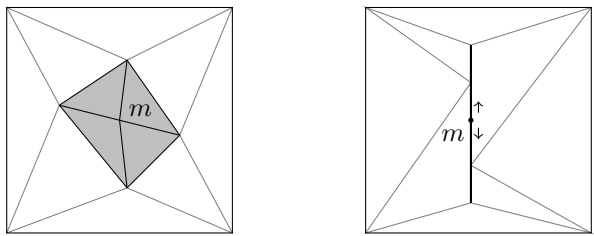
\begin{figure}[ht]
\begin{center}
\begin{tikzpicture}[-]

\matrix (m) [row sep = 8em, column sep = 5em]
{

\coordinate (P) at (0,0);
\coordinate (Q) at (3,0);
\coordinate (R) at (3,3);
\coordinate (S) at (0,3);
\coordinate (1) at (1.6,.6);
\coordinate (2) at (2.3,1.3);
\coordinate (3) at (1.6,2.3);
\coordinate (4) at (.7,1.7);
\coordinate (5) at (1.5,1.5);

\draw[gray] (P)--(1)--(Q)--(2)--(R)--(3)--(S)--(4)--cycle;
\draw (P)--(Q)--(R)--(S)--cycle;
\fill[gray!50!white] (1)--(2)--(5)--cycle;
\fill[gray!50!white] (3)--(2)--(5)--cycle;
\fill[gray!50!white] (3)--(4)--(5)--cycle;
\fill[gray!50!white] (1)--(4)--(5)--cycle;

\draw (1)--(2)--(3)--(4)--cycle;
\foreach \x in {1,2,3,4} \draw (5)--(\x);

\draw (5) node[anchor=205]{$m$};

&

\coordinate (1) at (1.4,.4);
\coordinate (2) at (1.4,.9);
\coordinate (3) at (1.4,2.5);
\coordinate (4) at (1.4,2);
\coordinate (5) at (1.4,1.5);

\draw[gray] (P)--(1)--(Q)--(2)--(R)--(3)--(S)--(4)--cycle;
\draw (P)--(Q)--(R)--(S)--cycle;

\draw[thick] (1) --(3);

\fill[black] (5) node[anchor=40]{$m$} circle (1pt);

\draw[->] (1.5,1.6) -- (1.5,1.75);
\draw[->] (1.5,1.4) -- (1.5,1.25);

\\};

\end{tikzpicture}
\caption{The vertex $m$ is a mosquito in the $(T,C)$ shown on the left.}
\label{fig:mosquito}
\end{center}
\end{figure}

Second, we call a vertex $v\in C$ a \emph{mosquito} if all neighbors of $v$
are also in $C$.  (In particular $v$ is not in $\bdry T$.)  If $v$ is a
mosquito then we may
form a new triangulation $\T'=(T',C')$ as follows:  $T'$ is obtained 
from $T$ by deleting $v$ and re-triangulating
the star of $v$ with no interior vertices, and $C'$ is obtained from 
$C$ by deleting $v$.  Let $f'_{\pentagon}$ be the corresponding
area map for $\T'$, and note that the image of $f'_{\pentagon}$ is exactly equal
to the image of $f_{\pentagon}$.  Thus the presence of $v$ in $T$ gives rise
to an independent dimension in the fiber of $f_{\pentagon}$.  See Figure \ref{fig:mosquito}.

\begin{thm}\label{thm:dimension}
Let $\T=(T,C)$ be a triangulation with $|\bdry T|=n$,
and let ${\pentagon}$ be an $n$-gon shape.
Assume $C$ is non-separating, and let $m$ be the number of mosquitos in $\T$.  
Then:
\begin{enumerate}
\item[(a)] 
the dimension of $V_{\pentagon}(\T)$ is $\dim X_{\pentagon} - m - 6$;
\item[(b)] 
the codimension of $V_{\pentagon}(\T)$ in $Y$ is $n-3$.
\end{enumerate}
\end{thm}

We now prove some preliminary lemmas. 

\begin{lem}
Let $\T=(T,C)$ be a triangulation.  If $C$ is non-separating then $\Triangles(C)$
is contiguous, i.e., the graph $C^*$ is connected.  If in addition $\T$ has no 
mosquitos then $C^*$ is a tree.
\end{lem}

\begin{proof}
Assume $C$ is non-separating, and suppose $C^*$ is not connected.  
Let $S$ be a contiguous set of triangles such that $C=\Vxs(S)$.
The union of the triangles of $\Triangles(C)$ is a connected 
subset of $T$.  Since $S$ is contiguous, $S^*$ is connected and therefore contained
in a single component of $C^*$, called the main component.  

Let $A$ be the union of the triangles making up a component of $C^*$
other than the main component.
For each edge $e$ of $\bdry A$ that is not in $\bdry T$, there is a vertex $z_e$
of $T$ such that $\{z_e\}\cup\bdry e$ is a triangle of $T\setminus A$.  See Figure \ref{fig:nonnonsep}.
Note that
$z_e$ is not in $C$.  Let $e=xy$ be one such edge.  Note both $x$ and $y$ are 
in $C$, and neither of the triangles containing the edge $xy$ is in $S$.  
Thus there is an arc from $x$ to $y$
in $S$, and this arc together with the edge $xy$ forms a loop 
$\gamma$ in $T$.  This loop encloses either $z_e$ or all the $z_f$ for $f\ne e$.
In either case the hypothesis that $C$ is non-separating is contradicted.
We conclude that $C^*$ is connected.

\begin{figure}[ht]
\begin{center}
\begin{tikzpicture}[scale=2, rotate=90]

\coordinate (x) at (2.2,1.7);
\coordinate (y) at (2.1,1.1);
\coordinate (z) at (2.6,1.4);
\coordinate (v) at (1.6,1.5);
\coordinate (1) at (2.2,2.3);
\coordinate (2) at (1.9,2.5);
\coordinate (3) at (1.8,2.2);
\coordinate (4) at (1.5,2.4);
\coordinate (5) at (1.6,2.1);
\coordinate (6) at (1.2,2.2);
\coordinate (7) at (1.1,1.7);
\coordinate (8) at (1.2,1.3);
\coordinate (9) at (1.5,1);
\coordinate (10) at (1.2,0.8);
\coordinate (11) at (1.6,0.5);
\coordinate (12) at (1.4,0.1);
\coordinate (13) at (1.8,0);
\coordinate (14) at (2.1,0.2);
\coordinate (15) at (1.9,0.6);
\coordinate (16) at (2.3,0.7);
\coordinate (17) at (1.8,1);
\coordinate (18) at (1.8,1.9);

\fill[gray!50!white] (1)--(x)--(3)--cycle;
\fill[gray!50!white] (1)--(2)--(3)--cycle;
\fill[gray!50!white] (4)--(2)--(3)--cycle;
\fill[gray!50!white] (4)--(5)--(3)--cycle;
\fill[gray!50!white] (4)--(5)--(6)--cycle;
\fill[gray!50!white] (7)--(5)--(6)--cycle;
\fill[gray!50!white] (7)--(5)--(v)--cycle;
\fill[gray!50!white] (7)--(8)--(v)--cycle;
\fill[gray!50!white] (9)--(8)--(v)--cycle;
\fill[gray!50!white] (9)--(8)--(10)--cycle;
\fill[gray!50!white] (9)--(11)--(10)--cycle;
\fill[gray!50!white] (12)--(11)--(10)--cycle;
\fill[gray!50!white] (12)--(11)--(13)--cycle;
\fill[gray!50!white] (15)--(11)--(13)--cycle;
\fill[gray!50!white] (15)--(14)--(13)--cycle;
\fill[gray!50!white] (15)--(14)--(16)--cycle;
\fill[gray!50!white] (15)--(y)--(16)--cycle;

\draw (1)--(x) node[anchor=300]{$x$}--(3)--(1)--(2)--(3)--(4)--(2);
\draw (3)--(5)--(4)--(6)--(5)--(7)--(6);
\draw (5)--(v) node[anchor=35]{$v$}--(7)--(8)--(v)--(9)--(8)--(10)--(9)--(11)--(10)--(12)--(11)--(13)--(12);
\draw (11)--(15)--(13)--(14)--(15)--(16)--(14);
\draw (15)--(y) node[anchor=260]{$y$}--(16);

\draw[gray] (x)--(y)--(v)--cycle;

\foreach \p in {9,11,15,v,y} \draw[gray] (\p)--(17);
\foreach \p in {3,5,v,x} \draw[gray] (\p)--(18);

\draw[gray] (x)--(z) -- (y);
\draw (z) node[anchor=south]{$z_e$};

\end{tikzpicture}
\caption{If $S$ is the set of shaded triangles, then the condition $C=\Vxs(S)$ is not non-separating.
Note that $S$ is contiguous but $\Triangles(C)$ includes the triangle $vxy$ and is not contiguous.}
\label{fig:nonnonsep}
\end{center}
\end{figure}
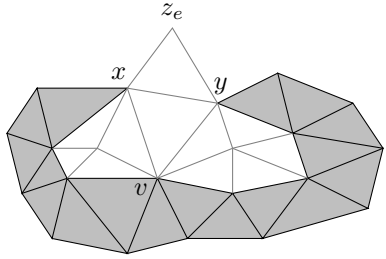


Finally, observe that any loop in $C^*$ would enclose some vertices of $T$.  
If all such vertices are in $C$ they would be mosquitos, but 
if some vertex is not in $C$ then $C$ would fail to be non-separating.  
Thus if $\T$ has no mosquitos then $C^*$ is a tree.
\end{proof}

\begin{lem}\label{lem:c=c^*}
Let $\T=(T,C)$ be a triangulation with no mosquitos and with $C$ non-empty and 
non-separating.
Then $|\Triangles(C)|=|C|-2$.
\end{lem}

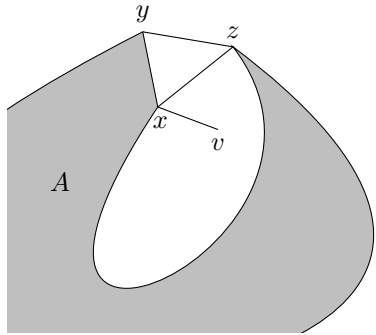
\begin{figure}[ht]
\begin{center}
\begin{tikzpicture}[-]

\coordinate (x) at (2,3);
\coordinate (y) at (1.8,4);
\coordinate (z) at (3,3.8);
\coordinate (v) at (2.8,2.7);

\clip (0,0) rectangle (6,4.5);

\draw  (y) node[anchor=south]{$y$} -- (z) node[anchor=south]{$z$}-- (x) node[anchor=100]{$x$};
\filldraw[fill=gray!50!white]
	(x) .. controls +(-3,-4.5) and (5,1.2) .. (z)
	.. controls +(8,-6) and (-10,-2) .. (y) -- cycle;

\draw (.7,2) node{$A$};
\draw (x) -- (v) node[anchor=north]{$v$};

\end{tikzpicture}
\caption{A new triangle $xyz$ could not already have $z\in C$.}
\label{fig:t+2=v}
\end{center}
\end{figure}

\begin{proof}
Note that $\Triangles(C)$ consists of a single triangle if and only if $|C|=3$.

By the previous lemma, the graph $C^*$ is a tree.  
We think of building the tree $C^*$ by starting with a vertex and adding leaves.
Suppose at some step we have a collection of triangles $A$, and we wish to
add a new triangle $xyz$.
The new triangle contributes 1 to $|\Triangles(C)|$.  The new triangles contributes at most 
1 to $|C|$, because the triangle $xyz$ shares an edge, say $xy$, with $A$.
We now show that $z$ was not already in $A$, and hence the triangle
$xyz$ contributes exactly 1 to $|C|$.

The argument is similar to the argument that $C^*$ contains no cycles; see Figure \ref{fig:t+2=v}.
Suppose on the contrary that $z$ is a vertex of a triangle of $A$.  Let $w$
be the midpoint of the edge $xy$.  Then $z$ is connected to $w$ by an arc
in $A$.  This arc together with the line segment from $w$ to $z$ in the triangle
$xyz$ makes a loop $\gamma$ that separates $x$ from $y$.  Relabel if
necessary so that $x$ is in the bounded region determined by $\gamma$.
Now, since $\T$ has no mosquitos, $x$ must have a neighbor $v$ that is not
in $C$.  But $v$ is also in the bounded region determined by $\gamma$,
contradicting the hypothesis that $C$ is non-separating.
\end{proof}

The proof of Theorem \ref{thm:dimension}
is an induction based on contracting an edge of $T$.
In order to define a contraction of a triangulation
$\T=(T,C)$, we need some language and notation.

First, we say that $T$ contains a \emph{subdivision} if there are vertices $x,y,z$ 
of $T$ such that $xy$, $yz$, $xz$ are edges of $T$, $xyz$ is not a triangle of $T$, 
and $\{x,y,z\}\ne \bdry\T$.  

Next, an edge $e=xy$ of $T$ is \emph{contractible} if $e\not\subset\bdry T$
and the ``contraction" $T/e$ (defined below) is a simplicial complex.  Specifically
this means that there are exactly two vertices $w,z$ of 
$T$ that are adjacent to both $x$ and $y$; when this holds, we define $T/e$
to be the simplicial complex obtained from $T$ by deleting the triangles $wxy$ and $xyz$
and identifying the vertices $x$ and $y$ 
(and consequently identifying the edges $wx$ and $wy$, and the edges $xz$ and $yz$). 

\begin{figure}[ht]
\begin{center}
\begin{tikzpicture}[-]

\matrix (m) [column sep = 5em]
{

\coordinate (x) at (2,4);
\coordinate (y) at (3,2);
\coordinate (z) at (1.2,2.4);
\coordinate (w) at (3.4,3.6);

\draw (x)--(z) node[anchor=0]{$z$}--(y)--(x)--(w)node[anchor=north west]{$w$}--(y);
\draw (2.2,3) node{$e$};

\draw (x)node[anchor=0]{$x$} -- +(0.2,0.6);
\draw (x) -- +(-0.4,0.4);

\draw (y)node[anchor=103]{$y$} -- +(0.6,0);
\draw (y) -- +(0.4,-0.4);
\draw (y) -- +(-0.2,-0.4);

&

\draw (x)--(z) node[anchor=0]{$z$}--(y)--(x)--(w)node[anchor=north west]{$w$}--(y);
\draw (2.2,3) node{$e$};

\draw (x)node[anchor=0]{$x$} -- +(0.2,0.6);
\draw (x) -- +(-0.4,0.4);

\draw (y)node[anchor=103]{$y$} -- +(0.6,0);
\draw (y) -- +(0.4,-0.4);
\draw (y) -- +(-0.2,-0.4);

\coordinate (v) at (5,4.5);
\draw (x)--(v)--(y);
\draw (v) node[anchor=west]{$v$};
\\};

\end{tikzpicture}
\caption{The figure on the left shows a contractible edge $e$, as long as there is no vertex $v$ 
as there is in the right figure.  
The edge $e$ on the right is not contractible, regardless of how the interior of the quadrilateral
$xwyv$ is triangulated.}
\label{fig:contractible}
\end{center}
\end{figure}
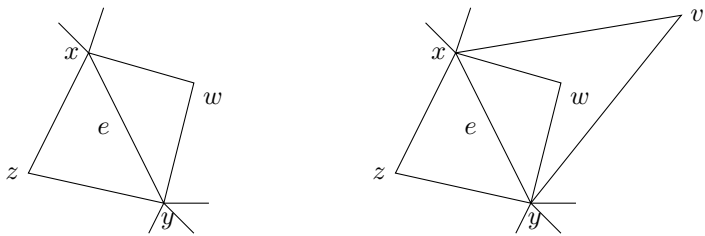

We leave it to the reader to verify that if the above condition holds, then $T/e$ is indeed simplicial.
See Figure \ref{fig:contractible}.

Note in the case that $e=xy$ is not a contractible edge, either $e$ is on the 
boundary or $T$ contains a subdivision.  Said differently, if $T$ has no subdivisions, then
every non-boundary edge of $T$ is contractible.

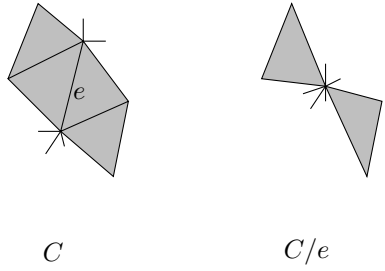
\begin{figure}[ht]
\begin{center}
\begin{tikzpicture}[-]

\coordinate (1) at (2.4,3.8);
\coordinate (2) at (2,2.8);
\coordinate (x) at (3,3.3);
\coordinate (y) at (2.7,2.1);
\coordinate (5) at (3.6,2.5);
\coordinate (6) at (3.4,1.5);

\matrix (m) [column sep = 5em]
{

\fill[gray!50!white] (2)--(1)--(x)--cycle;
\fill[gray!50!white] (2)--(y)--(x)--cycle;
\fill[gray!50!white] (y)--(5)--(x)--cycle;
\fill[gray!50!white] (y)--(5)--(6)--cycle;

\draw (2)--(1)--(x) --(2)--(y) -- (x)--(5)--(y)--(6)--(5);
\draw (x) -- +(0,0.3);
\draw (x) -- +(0.3,0);
\draw (y) -- +(-0.3,0);
\draw (y) -- +(.05,-0.2);
\draw (y) -- +(-0.2,-0.3);
\draw (2.95,2.6) node{$e$};
\draw (2.6,0.5) node{$C$};

&

\coordinate (xy) at (2.85,2.7);

\fill[gray!50!white] (2)--(1)--(xy)--cycle;
\fill[gray!50!white] (xy)--(5)--(6)--cycle;

\draw (2)--(1)--(xy)--cycle;
\draw (5)--(6)--(xy)--cycle;
\draw (xy) -- +(0,0.3);
\draw (xy) -- +(0.2,0.1);
\draw (xy) -- +(-0.3,-.1);
\draw (xy) -- +(0,-0.2);
\draw (xy) -- +(-0.2,-0.3);
\draw (2.6,0.5) node{$C/e$};

\\};
\end{tikzpicture}
\caption{If $C$ consists of the vertices of the four triangles shown on the left, then even 
if the edge $e$ is contractible, it is not good, because $C/e$ is not a condition.}
\label{fig:ungood}
\end{center}
\end{figure}

Suppose $\T=(T,C)$ is a triangulation and $e=xy$ is a contractible edge of $T$.
If $wxy$ and $xyz$ are triangles of $T$ (with $w\ne z$) then we define 
$C/e=\Vxs(\Triangles(C)\setminus\{wxy,xyz\})\subset\Vxs(T/e)$.
Thus $C/e=\emptyset$ if $C\subset\{w,x,y,z\}$, and otherwise $C/e$ is the image of $C$
under the quotient map $\Vxs(T)\to\Vxs(T/e)$.  Note however that $C/e$ may not 
be a condition since contracting $e$ may destroy contiguity.  See Figure \ref{fig:ungood}.

To ensure that $C/e$ is a condition, we introduce one more concept.

\begin{defn}[Good edge]\label{def:good}
Suppose that $\T=(T,C)$ is a triangulation with $C$ non-separating. 
An edge $e=xy$ of $T$ is \emph{good} if $e$ is contractible and either
\begin{description}
\item[Type 1]
$\{x,y\}\cap C=\emptyset$, or
\item[Type 2]
$e\subset \bdry C$, i.e., there exist $w,z\in\Vxs(T)$ adjacent to both $x$ and $y$, with  
$w\in C, z\notin C$.
\end{description}
\end{defn}

\begin{remark}
Note that if $T$ has at least one interior vertex and no subdivisions, then
$T$ has good edges.  The following lemma allows us to contract these edges.
\end{remark}

\begin{lem} \label{lem:goodisgood}
Let $\T=(T,C)$ be a triangulation with $C$ non-separating.
If $e$ is a good edge then 
$C/e\subset\Vxs(T/e)$ is a non-separating condition in the triangulation $T/e$.
\end{lem}

\begin{proof}
Let $S$ be a contiguous set of triangles of $T$ such that $C=\Vxs(S)$.

If $e$ is a type 1 good edge, then $C/e$ consists of the same vertices as $C$,
and $S$ remains contiguous in $T/e$.  Thus $C/e$ is a condition.
Also, the complement of $C/e$ in $T/e$ is just a contraction of the 
complement of $C$ in $T$.  Since $C$ is non-separating, so is $C/e$.

Now suppose $e$ is a type 2 good edge.  Let $S/e$ be the image of $S$ under
the contraction.  Then $S/e$ is contiguous, and $C/e=\Vxs(S/e)$.
So $C/e$ is a condition.  Also, the complement of $C/e$ in $T/e$ is 
identical to the complement of $C$ in $T$, so $C/e$ is non-separating.
\end{proof}

Thus the good edges are the ones we will contract:  if $\T=(T,C)$ and $e$ are
as in the lemma, we define $\T/e$ to be the triangulation $(T/e,C/e)$.
We also sometimes denote this by $\T|_{x=y}$.  

We are now ready to prove Theorem \ref{thm:dimension}.

\begin{proof}[Proof of Theorem \ref{thm:dimension}:]
First, we show that it is possible to eliminate mosquitos.  Suppose that 
$v\in\Vxs(T)$ is a mosquito.  Let $<$ be the total order on $\Vxs(T)$ used in the proof of
Theorem \ref{irr2}, and let $\a$ be defined as it was in that proof.  Without loss of generality, 
we can assume that $v$ is the last vertex with respect to the chosen ordering.  It follows  that $\a(v)=1.$  
As we mentioned in the discussion preceding the theorem, we may form a new 
triangulation by eliminating $v$ and retriangulating the star of $v$.  This reduces 
the dimension of $X_{\pentagon}$ by $\a(v)=1$, reduces $m$ by 1, and leaves 
$Y$ and the image of $f$ unchanged.  Thus, by induction, we may assume that 
$m=0$.

We claim that for any $\T$ satisfying the hypotheses, conclusions (a) and (b)
are equivalent.  Because $m=0$, this follows from the equation $\dim Y=\dim X+n-9$, 
which we now verify.
Suppose that $C=\emptyset$.  Then $\dim X=\sum \a(v)=6+2k$.
Also, recall that the number of triangles of $T$ is $2k+n-2$, so $\dim Y=2k+n-3$,
as desired.
On the other hand, if $C\ne\emptyset$, then $\dim X =\sum \a(v)=6+2k-(|C|-2)$.
Also, $\dim Y=2k+n-3-|\Triangles(C)|$.  The claim now follows from Lemma 
\ref{lem:c=c^*}.

We now prove (a) and (b) together by induction on the number $k$ of interior 
vertices.  If $k=0$, then $C=\emptyset$ and $m=0$, so 
$\dim X_{\pentagon}= 6, \dim Y=n-3$, and $f$ is a constant map.
So (a) and (b) hold.

Suppose that $\T$ contains a subdivision with vertices $x,y,z$.  We may assume that not
all of $x,y,z$ are in $C$, by the hypothesis that $C$ is non-separating.  In this case, we can
form new triangulations $\T'$ and $\T''$ as follows.
We get $T'$ by deleting all vertices (and triangles) in the interior
of the triangle $xyz$ and replacing them with the single triangle $xyz$;
note that $\bdry T'$ is an $n$-gon.  We let $C'=\Vxs(\Triangles(C)\cap\Triangles(T'))$
and $\T'=(T',C')$.

We get $T''$ by taking the sub-triangulation of $T$ consisting of those triangles
inside $xyz$.  The boundary of $T''$ is a 3-gon.  Let $C''=\Vxs(\Triangles(C)\cap\Triangles(T''))$
and $\T''=(T'',C'')$.

There is a projection $\hat\pi$ from $Y$ to $Y'$ that replaces the coordinates
corresponding to triangles of $T''$ with their sum.  Note that $\hat\pi$ restricts 
to a projection
$\pi:V\to V'$.  Let $v'$ be a generic point of $V'$, i.e., $v'\in\im (f'_{\pentagon})$ and the 
coordinate corresponding to the sum of triangles from $T''$ is nonzero.  (We can do
this because this coordinate is not always zero, which follows from the assumptions
that $m=0$ and $C$ is non-separating.)
We claim that $\pi^{-1}(v')=\hat\pi^{-1}(v')$.  Clearly $\pi^{-1}(v')\subseteq \hat\pi^{-1}(v')$.
Let $v\in\hat\pi^{-1}(v')$; we wish to show $v\in V$.

To see this, let $\rho'\in X'$ with 
$f'_{\pentagon}(\rho')=v'$.  In the drawing $\rho'$, the triangle $xyz$ has a shape $\nabla$.
Consider the map $f_\nabla$.  By induction the codimension of its image is zero, which
means there is a drawing of $\nabla$ realizing the areas $(v_\Delta)$ where $\Delta$
runs through the triangles of $T''$.  Inserting this drawing of $\nabla$ produces a drawing
$\rho\in X$ with $f_{\pentagon}(\rho)=v$.  

Now that we know that $\pi^{-1}(v')=\hat\pi^{-1}(v')$, we conclude that 
$\dim V=\dim V'+\dim F=\dim V' + \dim V''$.  Since we also have
$\dim Y=\dim Y'+\dim Y''$, we get $\codim_Y V = \codim_{Y'}V' + \codim_{Y''}V''
=\codim_{Y'}V'$.
Now the induction hypothesis applied to $\T'$ implies
that $\codim_Y V= n-3$, proving that (b) (and hence also (a)) holds for $\T$.

The remaining case is that $\T$ contains no subdivisions.
We will prove (a) (under the assumption $m=0$); we first note that the inequality 
$\dim  V_{\pentagon}(\T) \le  \dim X_{\pentagon} - 6$ follows from the previous 
discussion about the fibers of $f$.  

For the opposite inequality, recall that we are in the following situation.  The triangulation
$\T$ has no mosquitos, no subdivisions, and some number $k>0$ of interior vertices.
As we remarked earlier, in this setting there exists a good edge $e=xy$ of $T$.
By Lemma \ref{lem:goodisgood}, the induction hypothesis applies to the contraction $\T/e$,
which we now denote $\T'$.

Furthermore, from the definition of goodness it follows that $\T'$ has no mosquitos.  
Thus, inductively, the image of $f'_{\pentagon}$ has dimension at 
least $\dim X'-6$.  Moreover $\dim X=\dim X' +2$, since $y\not\in C$ so $\a(y)=2$.
Hence it suffices to show that $\dim \im( f_{\pentagon}) \geq \dim \im(f'_{\pentagon}) +2$.

Now let $w$ and $z$ be the vertices of $T$ that are adjacent to both $x$ and $y$.
Choose $\rho'\in X'$ such
that $w, x, z$ are not collinear and such that the rank of $Df'_{\pentagon}(\rho')$ is maximal. 
This can be done because by the definition of goodness, we may assume $w,x,z$ are
not all contained in $C$.
We view $X'$ as the subset of $X$ where $\rho(y)=\rho(x)$, and $f'_{\pentagon}$ as a restriction of
$f_{\pentagon}$.

Suppose $e$ is a type 1 good edge.  Then $y\notin C$ so $\a(y)=2$, and
we wish to show that the rank of $Df_{\pentagon}$ at $\rho'$ is at least two more than the 
rank of $Df'_{\pentagon}$ at $\rho'$.
This is achieved by moving $y$ off $x$ in the two independent directions $\rho(w)-\rho(x)$ and 
$\rho(z)-\rho(x)$.  These tangent directions in $X$ map via $Df_{\pentagon}$ to tangent
directions in $Y$ that are independent of each other and of the image of $Df'_{\pentagon}$, because
in the image of $Df'_{\pentagon}$, the areas of both $wxy$ and $xyz$ are constantly zero, and
moving in the direction $\rho(w)-\rho(x)$ keeps $wxy$ at zero but changes $xyz$ 
whereas moving in the direction $\rho(z)-\rho(x)$ keeps $xyz$ at zero but changes $wxy$. 
See Figure \ref{fig:alpha2}.
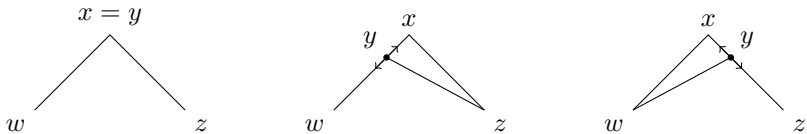
\begin{figure}[ht]
\begin{center}
\begin{tikzpicture}[-]

\matrix (m) [row sep = 8em, column sep = 3em]
{

\coordinate (x) at (0,2);
\coordinate (w) at (-1,1);
\coordinate (z) at (1,1);

\draw (w) node[anchor=north east]{$w$} -- (x) node[anchor=south]{$x=y$} -- (z) node[anchor=north west]{$z$};

&

\draw (w) node[anchor=north east]{$w$} -- (x) node[anchor=south]{$x$} -- (z) node[anchor=north west]{$z$};

\draw[<->] (-.45,1.55) --  (-.15,1.85);
\draw[fill=black] (-.3,1.7) node[anchor=south east]{$y$} circle (1pt);
\draw (z) -- (-.3,1.7);

&

\draw (w) node[anchor=north east]{$w$} -- (x) node[anchor=south]{$x$} -- (z) node[anchor=north west]{$z$};

\draw[<->] (.45,1.55) --  (.15,1.85);
\draw[fill=black] (.3,1.7) node[anchor=south west]{$y$} circle (1pt);

\draw (w) -- (.3,1.7);

\\};

\end{tikzpicture}
\caption{The case $\a(y)=2$:  as $y$ moves off of $x$ in the indicated directions, it generates two 
new tangent directions in area space.}
\label{fig:alpha2}
\end{center}
\end{figure}


On the other hand if $e$ is a type 2 good edge, then $x,y\in C$ and (after relabeling
if necessary) we may assume that $w\not\in C$ and
the order $<$ is such that $\a(y)=1$.  In this case we only need to show that the rank of 
$Df_{\pentagon}$ at $\rho'$ is at least one more than the 
rank of $Df'_{\pentagon}$ at $\rho'$.  Here, this is achieved by moving $y$ off $x$ in the direction of
the line $\ell$ associated to $C$.  This changes the area of $wxy$ and hence is independent
of all directions in the image of $Df'_{\pentagon}$.

This establishes the inequality $\dim f \geq \dim X -m-6$ and shows that (a) (and hence (b))
holds for $\T$.
\end{proof}

\section{Squares and the Monsky polynomial} \label{sec:Monsky}

We now turn our attention to triangulations of a square, i.e., $|\bdry T|=4$.
We rename the corners $P,Q,R,S$.
By Theorem \ref{thm:dimension}, if $\T=(T,C)$ is a triangulation with
$|\bdry T|=4$, then for any boundary shape, the (irreducible) projective variety 
$V_{\pentagon}$ has codimension 1 inside the projective space $Y$.

\begin{cor}\label{cor:existence}
Let $\T=(T,C)$ be a triangulation with $|\bdry T|=4$.  Then for any $4$-gon shape
$\pentagon$,
$V_{\pentagon}(\T)$ is the zero set of a single irreducible homogeneous 
polynomial.
\end{cor}

\begin{defn}[Monsky polynomial and degree]
Let $\T=(T,C)$ be a triangulation with $|\bdry T|=4$,
and let $\square$ be any nondegenerate parallelogram.  
We define the \emph{Monsky polynomial} $p(\T)$ to be the unique irreducible polynomial 
that vanishes on $V(\T):=V_{\square}(\T)$,
scaled so the coefficients are integers with no common factor.  This polynomial is well-defined 
up to sign.  We also define the \emph{degree}
of $\T$ to be $\deg(\T)=\deg(p(\T))$.
\end{defn}

Note that the scaling of $p(\T)$ in the above definition is possible since $V(\T)$ is the image
of the map $f:=f_\square(\T)$, whose coordinate functions are polynomials with rational
coefficients.

\section{A degree algorithm}\label{sec:degree}

In this section we study the Monsky polynomial $p(\T)$.  In particular, we describe an 
algorithm for computing a lower bound on the degree of $p(\T)$.
The degree $\deg(\T)$ is an integer invariant of $\T$.  At present we do not know how $\deg(\T)$
is related to any classical graph-theoretic invariants of the 1-skeleton of $T$.

The idea of the algorithm is to intersect $V$ with various coordinate hyperplanes,
and consider the components of the intersection inductively.  
We can identify all the components of intersection, and the sum of the degrees of
these components gives a lower bound for the degree of $V$.
If each component occurs with multiplicity one, then this sum is exactly the degree 
of $V$.  However, it is difficult to identify the multiplicities of the components, which is
why we can prove only that our 
algorithm produces a lower bound on the degree.
On the other hand, the algorithm involves a choice of ordering of the triangles of $T$,
and in all cases we have observed, there is an ordering that produces the
correct degree.  We conjecture that this is always the case.  
See the discussion at the end of the section.

For the rest of this section, we let $\T=(T,C)$ be given, where $|\bdry T|=4$ and $C$ is a 
non-separating condition.  Let ${\pentagon}$ be a square.  By Corollary \ref{cor:existence}, 
$V(\T)$ is the zero set of a single polynomial equation $p(\T)$ in $Y$; we denote its degree by
$\deg(\T)$.

Our algorithm is recursive, and we will need a way to detect triangulations with $\deg \T=1$. 
For this purpose, we introduce two types of triangulations where it is easy to see that the equation 
is linear.   Later, we shall see that any $\T$ for which $\deg \T = 1$ is of one of these types.  

For a triangulation $\T=(T,C)$, we say that $x,y\in\Vxs(T)$ are \emph{joined} if
either $\{x,y\}\subset C$ or there is an edge $xy$ in the 1-skeleton of $T$.

\begin{defn}[Types D and X]
Let $\T=(T,C)$ be a triangulation.  We say that $\T$ is \emph{linear of type D} if there are
two opposite corners of $T$ that are joined.  We say that $\T$ is \emph{linear of type X} if
there is a vertex of $T$ that is joined to all four corners of $T$.
\end{defn}

If $\T$ is linear of type D or X, then there is an obvious linear relation among the areas, so 
$\deg \T = 1$.  (The letter D stands for ``diagonal''; for an example of type X see Figure \ref{fig:X}.)

The next lemma shows that we can eliminate subdivisions without affecting the degree.
Suppose $T$ contains a subdivision with vertices $x,y,z$; then as in the proof of 
Theorem \ref{thm:dimension} we can form a new triangulation $\T'=(T',C')$ as follows.  
We get $T'$ by deleting all vertices (and triangles) in the interior of the triangle $xyz$
and replacing them with the single triangle $xyz$, and we set $C'=C\cap\Vxs(T')$.

\begin{lem}\label{lem:subdivisions}
Let $\T=(T,C)$ be a triangulation with $C$ non-separating.  Suppose $T$ contains
a subdivision with vertices $x,y,z$ and define $\T'$ as above.  
Then the polynomial $p(\T)$ is obtained from $p(\T')$ by replacing the variable
$A_{xyz}$ with the linear form $\sum \Delta_i$ where the sum is over the triangles
in $T$ interior to $xyz$.  In particular, $\deg(\T)=\deg(\T')$.
\end{lem}

\begin{proof}
As in the proof of Theorem \ref{thm:dimension}, we have a projection $\hat\pi:Y\to Y'$ that restricts
to a projection $\pi:V\to V'$.  Note that $p(\T')\circ\hat\pi$ defines an irreducible 
homogeneous polynomial on $Y$.  This polynomial vanishes on $V$ because
$v\in V$ implies $\hat\pi(v)=\pi(v)\in V'$.  Therefore this polynomial equals $p(\T)$.
\end{proof}

A key step in our algorithm involves a choice of triangle to kill.  The next
lemma asserts that in most situations we can choose such a triangle, preserving the
non-separating property of the condition.

\begin{defn}[Killable triangle]\label{def:killable}
Let $\T=(T,C)$ be a triangulation with $C$ non-separating.  A triangle $\Delta$ of $T$ 
is called \emph{killable} if $\Delta\not\in\Triangles(C)$ and $C\cup\Vxs(\Delta)$ is a 
non-separating condition.
\end{defn}

\begin{lem}\label{lem:killable}
Fix $\T=(T,C)$ with $C$ non-separating and such that $T$ has no subdivisions.
Then there exists a killable triangle in $\T$, or else $\T$ is linear of type D or X. 
\end{lem}

\begin{figure}[ht]
\begin{center}
\begin{tikzpicture}[-]

\coordinate (P) at (0,0);
\coordinate (Q) at (3,0);
\coordinate (R) at (3,3);
\coordinate (S) at (0,3);
\coordinate (1) at (.3,.7);
\coordinate (2) at (.7,1.2);
\coordinate (3) at (1.2,1.6);
\coordinate (4) at (1.8,1.6);
\coordinate (5) at (2.3,1.2);
\coordinate (6) at (2.7,.7);

\fill[gray!50!white] (P) -- (1) -- (2) -- (3) -- (4) -- (5) -- (6) -- (Q) --cycle;
\draw (P) -- (1) -- (2) -- (3) -- (4) -- (5) -- (6) -- (Q) ;

\foreach \x in {1,2,3} \draw (\x) -- (S);
\foreach \x in {3,4,5,6} \draw (\x) -- (R);

\draw (P) -- (Q) -- (R) -- (S) -- cycle;

\draw (1.5,.6) node{$C$};
\draw (P) node[anchor=north east]{$P$};
\draw (Q) node[anchor=north west]{$Q$};
\draw (R) node[anchor=south west]{$R$};
\draw (S) node[anchor=south east]{$S$};

\end{tikzpicture}
\caption{Sunshine.}
\label{fig:sunshine}
\end{center}
\end{figure}
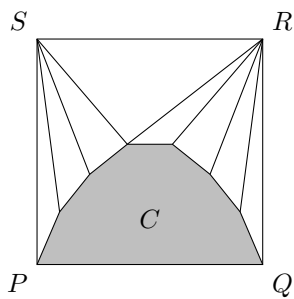

\begin{proof}
We will assume that $\T$ is not linear of type D or X, and we will show that there is a killable triangle.

If $C=\emptyset$ then any triangle that does not contain three corners of the
square is killable.  If each triangle contains three corners of the square then
the triangulation has exactly two triangles and $\T$ is linear of type D.  Henceforth, we assume 
$C\ne\emptyset$.

Recall that because $T$ has no subdivisions, every edge $e$ of $\bdry C$ that is not
contained in $\bdry T$ is good (see Definition \ref{def:good}).  For each such edge, there is a triangle 
$\Delta_e$ that is adjacent to a triangle of $\Triangles(C)$.  We will show that at least one of these
$\Delta_e$'s is killable.  Denote by $v_e$ the vertex of $\Delta_e$ that is not in $C$;
then $\Vxs(\Triangles(C) \cup \Delta_e)=C\cup\{v_e\}$.

Note that $\Triangles(C) \cup \Delta_e$ is contiguous, for each $e$.  Thus
$\Delta_e$ fails to be killable only if either $C\cup\{v_e\}$ contains
at least three corners (and thus isn't a condition) or $C\cup\{v_e\}$ is separating.

Fix $e$ and suppose $C'=C\cup\{v_e\}$ is separating.  Then there is
a loop $\gamma$ in the 1-skeleton of $T$ that uses only vertices of $C'$ and that 
encloses at least one vertex $w\not\in C'$.  See Figure \ref{fig:hornofafrica}.  Now, note that any 
triangle $\Delta_f$ in the bounded region determined by $\gamma$ is killable (for $(T,C)$) 
unless $C\cup\{v_f\}$ is separating.  
By an innermost argument, not all such conditions can be separating.  Thus
$(T,C)$ has is a killable triangle.

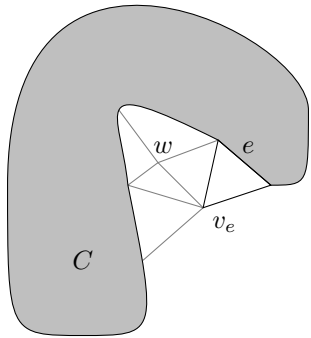
\begin{figure}[ht]
\begin{center}
\begin{tikzpicture}


\draw[gray] (1.6,2) -- (2.6,1.7) -- (1.8,1);
\draw[gray] (1.6,2) -- (2,2.3) -- (2.6,1.7);
\draw[gray] (2.8,2.6) -- (2,2.3) -- (1.48,3);

\filldraw[fill=gray!50!white] (0,2) .. controls (0,6) and (4,4) .. (4,3) 
	.. controls (4,2) and (4,2) .. (3.5,2) 
	-- (2.8,2.6)
	.. controls (1,3.5) and (1.5,3) .. (1.6,2)
	.. controls (2,0) and (2,0) .. (1,0)
	.. controls (0,0) and (0,0) .. (0,2);

\draw (3.5,2) -- (2.8,2.6) -- (2.6,1.7) -- cycle;

\draw (1,1) node{$C$};
\draw (3.2,2.5) node{$e$};
\draw (2.6,1.7) node[anchor=north west]{$v_e$};
\draw (2,2.3) node[anchor=250]{$w$};

\end{tikzpicture}
\caption{If $C\cup v_e$ is separating, then there is a killable triangle in the bounded region
separated by $C\cup v_e$.}
\label{fig:hornofafrica}
\end{center}
\end{figure}


We may therefore assume that for each $e$, $C\cup\{v_e\}$ is non-separating.
Now, if $C$ contains at most one corner, then $C\cup \{v_e\}$ contains at most
two corners and thus $\Delta_e$ is killable.  So, suppose $C$ contains two corners.  

If these corners are not adjacent, then $\T$ is 
linear of type D, so we may assume that the corners are adjacent, 
say $P$ and $Q$.
We need to show that for some $e\subset \bdry C$ other than $e=PQ$, 
$v_e$ is not a corner.  

But if $v_e$ is a corner for each such $e$, it follows that either all vertices of
$\bdry C$ are connected by an edge to $R$ (or $S$), in which case $T$ contains a
diagonal and is linear of type D, or else there is a vertex $v$ of $\bdry C$ that is 
connected by edges to both $R$ and $S$, in which case $\T$ is linear of type X.
See Figure \ref{fig:sunshine}.
\end{proof}


Armed with these two lemmas, we now present an algorithm that conjecturally
computes $\deg \T$.

\begin{algorithm}\label{alg}
Let the triangulation $\T=(T,C)$ be given, with $|\bdry T|=4$ and $C$ non-separating.

\begin{enumerate}
\item
If $\T$ has mosquitos, eliminate them by re-triangulating (as in the proof of
Theorem \ref{thm:dimension}).
\item
If $T$ has subdivisions, delete them to obtain $\T'$ as in Lemma \ref{lem:subdivisions}.  
Recursively set $d(\T)=d(\T')$.
\item
Look for a killable triangle.  If none exists, set $d(\T)=1$.  Otherwise, let $\Delta=xyz$ be a killable triangle.  If $C=\emptyset$, set $\T'=(T,\{x,y,z\})$ and recursively set $d(\T)=d(\T')$.
\item
If $C\ne\emptyset$, then $\Delta$ intersects $C$ in an edge $e$.  
Set $\T'=T/e$ and $\T''=(T,C\cup\{x,y,z\})$.
\item 
Discard the triangulation $\T''$ if it has any mosquitos. 
\item
If $\T''$ was not discarded in the previous step, recursively
set $d(\T)=d(\T')+d(\T'')$.  Otherwise, recursively set $d(\T)=d(\T')$.
\end{enumerate}
\end{algorithm}

\begin{remark}
Step (3) often allows for a choice of killable triangle $\Delta$.  Different choices
can lead to different outputs; see Section 6 for examples of this.  Thus $d(\T)$ is
not a well-defined integer.  However, we have the following result.
\end{remark}

\begin{thm}\label{thm:lowerbound}
Any output $d(\T)$ computed by Algorithm \ref{alg} satisfies $d(\T)\leq\deg(\T)$.
\end{thm}

\begin{proof}
We induct on the dimension of $V=V(\T)$.  Fix $\T=(T,C)$.  Eliminating mosquitos doesn't 
change $V$ (or the degree), so we may assume $\T$ has no mosquitos.  If $T$ has a subdivision then
we can delete it by Lemma \ref{lem:subdivisions}, and then by induction we are done.  
So we may also assume that $T$ has no subdivisions.  Next, if $\T$ has no killable triangles 
then the algorithm sets $d(\T)=1$, and by Lemma \ref{lem:killable} $\T$ is linear of type D or X,
so $\deg(\T)=1$ and again 
we are done.  We may therefore assume $\T$ has a killable triangle $\Delta=xyz$.

Consider the variety $\widetilde V=V\cap \{A_{\Delta}=0\}$.  If $C=\emptyset$ then 
$\widetilde V=V(\T')$ where $\T'=(T,\{x,y,z\})$, so $\widetilde V$ is irreducible by 
Theorem \ref{irr2}.
Thus the polynomial $p(\T)|_{\Delta=0}$ has the unique irreducible
factor $p(\T')$.  So $\deg p(\T) \ge \deg p(\T') \geq d(\T')=d(\T)$ (the 
last inequality by induction, the last equality by the algorithm).

If $C\ne\emptyset$, then there is an edge $xy$ in both $\Delta$ and $\bdry C$.  
In this case, a generic point
$q$ of $\widetilde V$ is in the image of the area map, and so there is a drawing
$\rho$ realizing the areas $q$.  Note that either $\rho(x)=\rho(y)$ or the four points
$\rho(x),\rho(y),\rho(w),\rho(z)$ are collinear (where $w$ and $z$ are the vertices
of $T$ that are adjacent to both $x$ and $y$).  Thus $q$ is in either $V'=V(\T')$ or $V''=V(\T'')$
(where we think of these as being subsets of $V$ in the obvious way).  Hence
the only possible components of $\widetilde V$ are $V'$ and $V''$, and each of
these is a component if and only if it has codimension 1 in $V$. 

To test whether $V'$ and $V''$ are components, we use Theorem \ref{thm:dimension}.  
Since $\dim X'=\dim X''=\dim X-1$, 
the only way that $V'$ or $V''$ could fail to have codimension 1 in $V$ is if the
corresponding triangulation has mosquitos.  Lemma \ref{lem:nomosquitos} (below) asserts
that $\T'$ has no mosquitos, so $V'$ is a component of $\widetilde V$.
If $\T''$ has mosquitos we discard it; otherwise $V''$ is also a component of $\widetilde V$.
In the case that both $V'$ and $V''$ are components of $\widetilde V$, we note that
$V'\ne V''$, since it is easy to construct points of one that are not on the other.
Therefore, using the induction hypothesis, if $\widetilde V$ has just one component, 
then $\deg p=\deg p|_{\Delta=0} \geq \deg \T'\geq d(\T')=d(\T)$, 
whereas if $\widetilde V$ has two components then 
$\deg p=\deg p|_{\Delta=0} \geq \deg p'+\deg p''\geq d(\T')+d(\T'')=d(\T)$.
\end{proof}

\begin{lem}\label{lem:nomosquitos}
Let $\T=(T,C)$ be a triangulation with no mosquitos and with $C$ nonempty and 
non-separating.  Let $\Delta=xyz$ be a killable triangle, where $\{x,y\}\subset C$.
Let $e=xy$ and let $\T'=T/e$.  Then $\T'$ has no mosquitos.
\end{lem}

\begin{proof}
If $C/e=\emptyset$ then certainly $\T'$ has no mosquitos.  Otherwise $C/e$ is
the image of $C$ under the identification $x=y$.  Suppose $v\in\Vxs(T/e)$ is a
mosquito in $\T'$.  Then $v$ and all its neighbors in $T/e$ are vertices of $C/e$.  
Since $x$ and $y$ are both in $C$, this implies that $v$ and all its neighbors in 
$T$ are vertices of $C$.  But $\T$ has no mosquitos, so this is a contradiction.
\end{proof}

\begin{remark}
In a recursive computation of $d(\T)$, if it happens that every time we intersect a
variety with a hyperplane, the components of intersection have multiplicity 1,
then each inequality in the above proof is an equality, and we have $d(\T)=\deg(\T)$.  
\end{remark}

\begin{conjecture}\label{conj:algworks}
For any $\T$, it is possible to run Algorithm \ref{alg} in such a way that
at each step, every component has geometric multiplicity one.
\end{conjecture}

If the conjecture is true, then by killing the triangles in the specified
order, the algorithm will produce the degree of $p(\T)$.

As an approach to this conjecture, we point out that if there are any nonsingular
points of $V$ in the intersection $V\cap \{A_\Delta=0\}$, then in most cases killing the triangle $\Delta$
results in components of multiplicity one.
To see this, one can use a simple differential geometry argument we call ``busting through.''  
The idea is to pick a nonsingular point $x$ of $V\cap\{A_\Delta=0\}$ and construct a curve
$\gamma\subset V$ through $x$ that is not tangent to $\{A_\Delta=0\}$.  The existence of 
$\gamma$ shows that $V$ itself is not tangent to $\{A_\Delta=0\}$, and so the multiplicity of 
intersection of any component must be one.  However, examples in Section \ref{sec:examples}
show that it is possible for $V\cap\{A_\Delta=0\}$ to consist entirely of singular points of $V$.

\section{Consequences of the degree algorithm}

Although we cannot guarantee that Algorithm \ref{alg} returns $\deg\T$, we are
nevertheless able to use the algorithm to deduce some results about the degree.
Specifically, in case $C=\emptyset$ we give a general lower bound for $\deg\T$ in 
terms of $|\Vxs(T)|$, and for arbitrary $C$ we characterize triangulations $\T$ with 
$\deg\T=1$.  The latter is similar to a result in \cite{gonewild}.

Our first corollary requires a graph-theoretic lemma.

\begin{lem}\label{lem:planargraphthing}
Let $T$ be a simplicial complex homeomorphic to a disk, with $|\bdry T|=4$ and
with at least one vertex in the interior.  Suppose $T$ has no subdivisions.  Then there
is an edge $e=xy$ of $T$ such that $\{x,y\}\not\subset\bdry T$ and $T/e$
has no subdivisions.
\end{lem}

\begin{proof}
We call a 4-cycle $\gamma$ in the 1-skeleton of $T$ \emph{bad} if $\gamma$ encloses 
at least one vertex of $T$ (equivalently, if $\gamma$ encloses more than 2 triangles).
Note that if the 4-cycle $wxyz$ is bad then neither $wy$ nor $xz$ is an edge of
$T$, since $T$ is assumed to have no subdivisions.

Let $e=xy$ be a non-boundary edge.  If $T/e$ contains a subdivision, then there must 
be a bad 4-cycle $\gamma$ including the edge $e$.  So it suffices to show that there is
an edge $e$ that is not part of any bad 4-cycle.  See Figure \ref{fig:bad4cycle}.

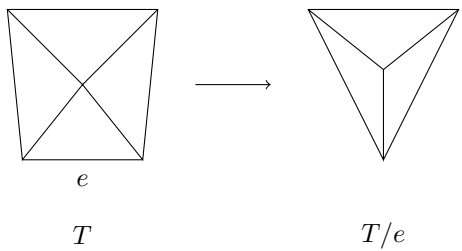
\begin{figure}[ht]
\begin{center}
\begin{tikzpicture}[-]

\draw (0,2) -- (1,1) --(.2,0) -- (0,2) -- (2,2) -- (1.8,0) -- (1,1) -- (2,2);
\draw (.2,0) -- node[below=2pt]{$e$} (1.8,0);
\draw (1,-1) node {$T$};

\draw [->] (2.5,1)--(3.5,1);

\draw (4,2) -- (5,1.2) --(5,0) -- (4,2) -- (6,2) -- (5,0);
\draw (6,2)--(5,1.2);
\draw (5,-1) node {$T/e$};

\end{tikzpicture}
\caption{The edge $e$ is part of a bad 4-cycle.}
\label{fig:bad4cycle}
\end{center}
\end{figure}


For bad 4-cycles $\gamma,\gamma'$ we say $\gamma\le\gamma'$ if every triangle
enclosed by $\gamma$ is also enclosed by $\gamma'$.
Fix a bad 4-cycle $\gamma=wxyz$ that is minimal with respect to this partial order.
Let $v$ be a vertex enclosed by $\gamma$, and let $e_0$ be an edge emanating from
$v$.  

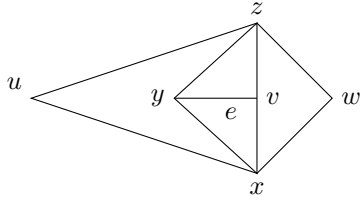
\begin{figure}[ht]
\begin{center}
\begin{tikzpicture}[rotate=90]

\coordinate (x) at (-1,0);
\coordinate (y) at (0,1.1);
\coordinate (u) at (0,3);
\coordinate (v) at (0,0);
\coordinate (z) at (1,0);
\coordinate (w) at (0,-1);

\draw (x) node[anchor=north]{$x$}
	-- (u) node[anchor=south east]{$u$}
	-- (z) node[anchor=south]{$z$}
	-- (w) node[anchor=west]{$w$}
	-- cycle;
\draw (x) -- (y) node[anchor=east]{$y$}
	-- (z) -- (v) node[anchor=west]{$v$}
	-- cycle;
\draw (v) -- node[anchor=north west]{$e$} (y);

\end{tikzpicture}
\caption{The edge $e$ is not part of a bad 4-cycle.}
\label{fig:cheating}
\end{center}
\end{figure}


If $e_0$ is not part of a bad 4-cycle, then we are done, so assume $e_0$ is part of a
bad 4-cycle $\beta$.  Since $\gamma$ is minimal, $\beta$ must enclose a triangle
not enclosed by $\gamma$.  Thus $\beta$ must intersect $\gamma$, and the four vertices of 
$\beta$ must be $v$, two of $w,x,y,z$, and a vertex $u$ that is outside $\gamma$.
Since $\beta$ is bad, neither of its diagonals is an edge of $T$, so the set of vertices
in $\beta\cap\gamma$ must be either $\{w,y\}$ or $\{x,z\}$.  Without loss of generality,
suppose $\beta=vxuz$.  Moreover, exactly one of $w,y$ is inside $\beta$; suppose it 
is $y$.  See Figure \ref{fig:cheating}.

Now, consider the 4-cycle $\alpha=vxyz$.  By minimality of $\gamma$, $\alpha$
does not enclose a vertex.  So $\alpha$ encloses exactly two triangles, and it follows
that one of the diagonals of $\alpha$ is an edge of $T$.  But $xz$ is not an edge
of $T$, because $\gamma$ is bad.  So $vy$ is an edge of $T$.  We will show that
$e=vy$ is not part of any bad 4-cycle.

Suppose there is a bad 4-cycle $\delta=vyab$.  It is impossible for $a$ to coincide
with $x$, because the edge $vx$ would be a diagonal of the bad 4-cycle
$vyxb$.  Likewise $a\ne z$ and $x\ne b \ne z$.  Thus, since $v$ is enclosed by 
$\gamma$ and $y$ is enclosed by $\beta$, the only possibility is $a=u$ and $b=w$.  

Now, the 4-cycle $\delta=vyuw$ separates $x$ from $z$ in $T$, and thus must 
enclose exactly one of these vertices.  But if $\delta$ encloses $x$, then $T$
contains the subdivision $wzu$, and if $\delta$ encloses $z$ then $T$ contains
the subdivision $wxu$.  This contradiction shows that $e$ satisfies the 
conclusion of the lemma.
\end{proof}

\begin{cor}\label{cor:linearlowerbound}
Let $\T=(T,\emptyset)$ be a triangulation of a square with $k$ interior vertices and
no subdivisions.  Then $\deg \T\ge k$.
\end{cor}

Example \ref{ex:diagonal} in the next section shows that this bound is sharp
for all $k$.

\begin{proof}
We proceed by induction; the conclusion is true if $k\le1$.  Suppose $k\ge 2$.

By Lemma \ref{lem:planargraphthing}, there is a non-boundary edge $e=xy$ such that $T/e$ has no 
subdivisions.  As $C=\emptyset$, the edge $e$ is good.  We run Algorithm \ref{alg}
as follows.  Let $w,z$ be the two vertices that are adjacent to both $x$ and $y$.
Observe that the triangle $wxy$ is killable.  We first kill $wxy$
and obtain $\T'=(T,\{w,x,y\})$ with $d(\T')=d(\T)$.  In $\T'$, the triangle $xyz$ is
killable, so we kill it, giving triangulations $\T''=\T/e$ and $\T'''=(T,\{w,x,y,z\})$.

Now, $\T''$ has $k-1$ interior vertices, no condition (hence no mosquitos), and by our 
choice of $e$, no subdivisions.  Thus by induction we have $\deg(\T'')\ge k-1$.
Also, since $T$ has no subdivisions, the condition in $\T'''$ is not big enough to 
create a mosquito, and so $d(\T''')\ge 1$.  Thus $\deg\T\ge \deg(\T'')+\deg(\T''')\ge k$.
\end{proof}

We need one more lemma before we characterize those triangulations $\T$ with $\deg \T=1$.

\begin{lem}\label{lem:xory}
Let $\T=(T,C)$ be a triangulation with no mosquitos and with $C$ non-separating.  
Suppose $\Delta=xyz$ is a killable triangle with $\{x,y\}\subset C$.  
Let $\T''=(T,C\cup\{z\})$.  If $\T''$ has a mosquito then either $x$ or $y$ is a mosquito of $\T''$.
\end{lem}

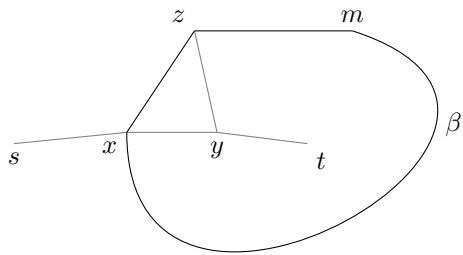
\begin{figure}[ht]
\begin{center}
\begin{tikzpicture}[scale=1.5]

\clip (-1,-1.5) rectangle (5,1.5);

\coordinate (x) at (.5,0);
\coordinate (y) at (1.3,0);
\coordinate (s) at (-.5,-0.1);
\coordinate (t) at (2.1,-0.1);
\coordinate (m) at (2.5,.9);
\coordinate (z) at (1.1,.9);

\draw[gray] (s) node[anchor=north,black]{$s$} 
	-- (x) node[anchor=north east,black]{$x$}
	-- (y) node[anchor=north,black]{$y$}
	-- (t) node[anchor=north west,black]{$t$};
\draw[gray] (z) -- (y);
\draw (x) .. controls (0.5,-2.5) and (5,.1) .. node[near end, anchor=west]{$\beta$} (m);
\draw (x) -- (z) node[anchor=south east]{$z$};
\draw (z) -- (m) node[anchor=south]{$m$};


\end{tikzpicture}
\caption{The loop $\beta$ must enclose either $s$ or $t$.}
\label{fig:xory}
\end{center}
\end{figure}


\begin{proof}
Let $m$ be a mosquito of $\T''$; we may assume $m\not\in\{x,y\}$.  
Thus $m\in C$ and since $\T$ has no mosquitos, there must
be an edge $mz$ in $T$.  Let $\alpha$ be an embedded path in $\bdry C$ from $m$ to $x$ such
that $\alpha$ does not go through $y$, and let $\beta$ be the union of $\alpha$ and the path $mzx$.  
Thus $\beta$ is an embedded loop in $T$.

Suppose neither $x$ nor $y$ is a mosquito in $\T''$.  Then there are neighbors $s,t$ of $x,y$ 
such that $s\ne z\ne t$ and $s,t\not\in C$.  See Figure \ref{fig:xory}. Consider the path $\gamma=sxyt$.  
Note that $\gamma$ intersects
$\beta$ in exactly one point, namely $x$.  This intersection is transverse, since the edges coming
out of $x$ occur in the cyclic order $y,z,s,\alpha$.  Thus exactly one of $s,t$ is enclosed by $\beta$.
But this means that $C\cup\{z\}$ fails to be non-separating, contradicting the assumption 
that $\Delta$ is killable.
\end{proof}

\begin{cor}\label{cor:linear}
Let $\T=(T,C)$ be a triangulation of a square with $C$ is non-separating.  Then 
$\deg\T=1$ if and only if $\T$ is linear of type D or X.
\end{cor}

\begin{proof}
We have already seen that triangulations that are linear of type D or X
have degree 1.  To prove the converse we will use induction on the dimension of $V(\T)$.  

Let $\T=(T,C)$ be a triangulation and suppose $\deg \T=1$.  We wish to show
that $\T$ is linear of type D or X.  If not, then by Lemma \ref{lem:killable}, when
we run Algorithm \ref{alg} we will find a killable triangle $\Delta=xyz$.

Consider first the case that $C=\emptyset$.  Then $\deg\T'=deg\T=1$ where
$\T'=(T,\{x,y,z\})$.  Since $\dim V'<\dim V$, it follows by induction that $\T'$
is linear of type D or X.  But two vertices are joined in $\T'$ if and only if they
are joined in $\T$; thus $\T$ is also linear of type D or X.  (Note that the case
$C=\emptyset$ also follows from Corollary \ref{cor:linearlowerbound}.)

Now suppose that $C\neq\emptyset$.  Let $\T'$ and $\T''$ be as in Algorithm \ref{alg}.
It follows from the proof of Theorem \ref{thm:lowerbound} that $\deg \T'=1$ and that
$\T''$ contains at least one mosquito.  By induction, $\T'$ is linear of type D or X.

We first consider the case that $\T'$ is linear of type D.  Thus two opposite corners,
say $P$ and $R$, are joined in $\T'$.  If $P$ and $R$ are joined in $\T$, we are done,
so assume the opposite.  Then it must be the case that either $x$ or $y$ is a corner,
say $x=R$, so that $y$ is connected to both $P$ and $R$.  Now, since $\T'$ has at
least one mosquito, Lemma \ref{lem:xory} implies that $y$ must be a mosquito in $\T''$.
In particular, this implies $P,R\in C$ and so $P$ and $R$ are joined in $\T$ after all.

Next suppose $\T'$ is linear of type X.  Let $v$ be a vertex of $\T'$ that is joined to
all four corners.  If $v$ is joined to all four corners in $\T$, then we are done, so assume
the opposite.  It follows that $v$ must be the image of $x$ (and $y$) under the contraction
$\Vxs(T)\to\Vxs(T/e)$.  Again, since $\T''$ has at least one mosquito, Lemma \ref{lem:xory}
implies that either $x$ or $y$ is a mosquito of $\T''$; suppose it is $x$.  Therefore every
neighbor of $x$ other than $z$ is in $C$.  We now claim that $y$ is joined to all four corners 
in $\T$, and hence $\T$ is linear of type X.  To see this, consider any corner $w$.  If $w=z$ then there is an edge from $w$ to $y$.  
Thus we can assume $w\ne z$.  Since $w$ is joined to $v$ in $\T'$, it must be that $w$ is 
joined to either $x$ or $y$ in $\T$.  But if $w$ is joined to $x$ in $\T$, then $w\in C$, because 
$w\ne z$ and every neighbor of $x$ in $\T$ other than $z$ is in $C$.  Thus either way, $w$ 
is joined to $y$ in $\T$.
\end{proof}

\section{Examples}\label{sec:examples}

\begin{example}\label{example2}
Let $\T=\T_2=(T_2,\emptyset)$ be the triangulation in Figure \ref{fig:example2}.  
Using elementary algebra one can easily compute that 
$p(\T)=(A+C+E)^2-4AC-(B+D+F)^2+4DF$; see Example \ref{ex:diagonal}
for a more general computation.  In particular $\deg(\T)=2$.
One way to run the algorithm on $\T$ is the following.  Note that $T$ has
no subdivisions.  We first kill $B$ to obtain $\T^1=(T,\{P,x,y\})$ (where $x$ and $y$ are the
names of the two interior vertices), and we set $d(\T)=d(\T^1)$.

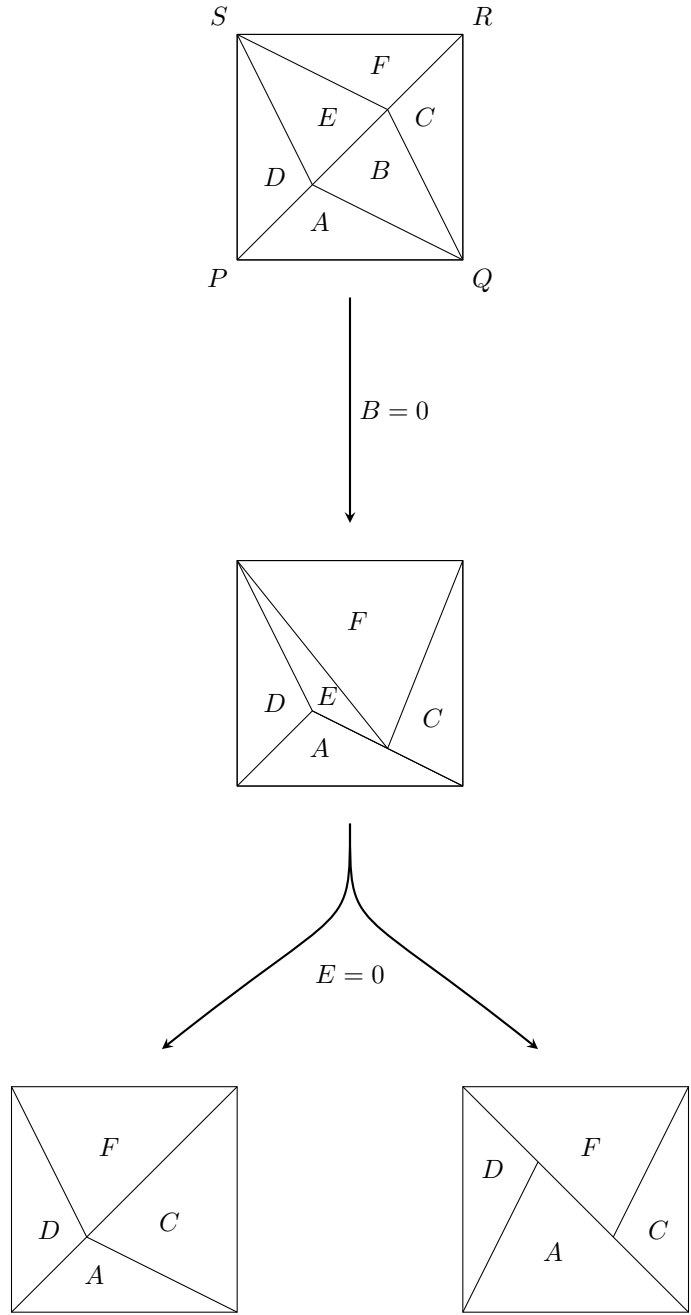
\begin{figure}[ht]
\begin{center}
\begin{tikzpicture}[>=stealth]

\coordinate (P) at (0,0);
\coordinate (Q) at (3,0);
\coordinate (R) at (3,3);
\coordinate (S) at (0,3);

\coordinate (x) at (1,1);
\coordinate (y) at (2,2);

\coordinate (A) at (0,14);

\draw ($(P) + (A)$) node[anchor=north east] {$P$} 
	-- ($(S) + (A)$) node[anchor=south east] {$S$}
	-- ($(R) + (A)$) node[anchor=south west] {$R$}
	-- ($(Q) + (A)$) node[anchor=north west] {$Q$}
	-- cycle;

\foreach \p in {P,x,y,R} \draw ($(S)+(A)$) -- ($(\p)+(A)$);
\foreach \p in {P,x,y,R} \draw ($(Q)+(A)$) -- ($(\p)+(A)$);
\draw ($(P)+(A)$)--($(x)+(A)$)--($(y)+(A)$)--($(R)+(A)$);
\draw ($(1.1,0.5) + (A)$) node{$A$};
\draw ($(1.9,1.2) + (A)$) node{$B$};
\draw ($(2.5,1.9) + (A)$) node{$C$};
\draw ($(0.5,1.1) + (A)$) node{$D$};
\draw ($(1.2,1.9) + (A)$) node{$E$};
\draw ($(1.9,2.6) + (A)$) node{$F$};


\coordinate (B) at (0,7);

\coordinate (y') at (2,0.5);

\draw ($(Q) + (B)$) -- ($(R) + (B)$) -- ($(S) + (B)$) -- ($(P) + (B)$) -- cycle;
	
\foreach \p in {P,x,y',R} \draw ($(S)+(B)$) -- ($(\p)+(B)$);
\foreach \p in {P,x,y',R} \draw ($(Q)+(B)$) -- ($(\p)+(B)$);
\draw ($(P)+(B)$)--($(x)+(B)$)--($(y')+(B)$)--($(R)+(B)$);

\draw ($(1.1,0.5) + (B)$) node{$A$};
\draw ($(2.6,0.9) + (B)$) node{$C$};
\draw ($(0.5,1.1) + (B)$) node{$D$};
\draw ($(1.2,1.2) + (B)$) node{$E$};
\draw ($(1.6,2.2) + (B)$) node{$F$};


\coordinate (C) at (-3,0);

\draw ($(Q) + (C)$) -- ($(R) + (C)$) -- ($(S) + (C)$) -- ($(P) + (C)$) -- cycle;
\draw ($(Q) + (C)$) -- ($(x) + (C)$) -- ($(S) + (C)$);
\draw ($(R) + (C)$) -- ($(x) + (C)$) -- ($(P) + (C)$);

\draw ($(1.1,0.5) + (C)$) node{$A$};
\draw ($(2.1,1.2) + (C)$) node{$C$};
\draw ($(0.5,1.1) + (C)$) node{$D$};
\draw ($(1.3,2.2) + (C)$) node{$F$}; 


\coordinate (D) at (3,0);
\coordinate (x'') at (1,2);
\coordinate (y'') at (2,1);

\draw ($(Q) + (D)$) -- ($(R) + (D)$) -- ($(S) + (D)$) -- ($(P) + (D)$) -- cycle;

\draw ($(S)+(D)$) -- ($(Q)+(D)$);
\draw ($(P)+(D)$) -- ($(x'')+(D)$);
\draw ($(y'')+(D)$) -- ($(R)+(D)$);

\draw ($(1.2,0.8) + (D)$) node{$A$};
\draw ($(2.6,1.1) + (D)$) node{$C$};
\draw ($(1.7,2.2) + (D)$) node{$F$};
\draw ($(0.4,1.9) + (D)$) node{$D$};


\draw[->,thick] (1.5,13.5) -- (1.5, 10.5);

\draw[->,thick] (1.5,6.5) .. controls (1.5,5) and (1.5,5.5) .. (-1,3.5);
\draw[->,thick] (1.5,6.5) .. controls (1.5,5) and (1.5,5.5) .. (4,3.5);


\draw (1.5,12) node[right] {$B=0$};
\draw (1.5,4.5) node{$E=0$};

\end{tikzpicture}
\caption{The top triangulation, $\T_2$, has degree 2.  Its Monsky polynomial 
is $(A+C+E)^2-4AC-(B+D+F)^2+4DF$.}
\label{fig:example2}
\end{center}
\end{figure}

Next kill $E$, obtaining $\T^2=(T|_{x=y},\emptyset)$ and $\T^3=(T,\{P,x,y,R\})$.
Each of these has the correct dimension, so we set $d(\T^1)=d(\T^2)+d(\T^3)$.
After killing one more triangle of $\T^2$, there will be no more killable triangles,
and the algorithm will return $d(\T^2)=1$, which in fact is the correct degree since
$p(\T^2)=A-C-D+F$.
Meanwhile, $\T^3$ already has no killable triangles, so the algorithm returns
$d(\T^3)=1$, and again this is the correct degree since $p(\T^3)=A-C+D-F$.

Thus the result of the algorithm is that $d(\T)=d(\T^1)=d(\T^2)+d(\T^3)=1+1=2$, 
and in this case the algorithm produces the correct answer, $d(\T)=\deg(\T)=2$.

Note that $p(\T^1)=p(\T)|_{B=0}$ is irreducible; this is why $\deg(\T^1)=\deg(\T)$.
On the other hand, $p(\T^1)|_{E=0}=(A+C)^2-4AC-(D+F)^2+4DF$ factors as
$(A-C-D+F)(A-C+D-F)$.  These factors are the polynomials associated to
the combinatorial ``factors'' $\T^2$ and $\T^3$ of $\T^1|_{B=0}$.
\end{example}

Up to isomorphism, $T_2$ is the only triangulation with two
interior vertices and no subdivisions.  Thus, by Corollary \ref{cor:linearlowerbound},
$\T_2$ is the only triangulation with $C=\emptyset$ that satisfies $\deg\T=2$.

\begin{example}[The diagonal case]\label{ex:diagonal}
Examples \ref{example1} and \ref{example2} begin an infinite sequence of
triangulations $\T_n=(T_n,\emptyset)$ that we refer to as the ``diagonal case.''  For $n\ge 1$, the 
triangulation $T_n$ has a total of $n+4$ vertices, with the $n$ interior vertices arranged
along a diagonal of the square.  Each is connected to its neighbors and to the two opposite
corners of the square, as in Figure \ref{fig:diagonal}.  (In viewing this figure, note that for $1\le i \le n$ 
the $p_i$ are not constrained to be on the diagonal of the square.)

\begin{center}
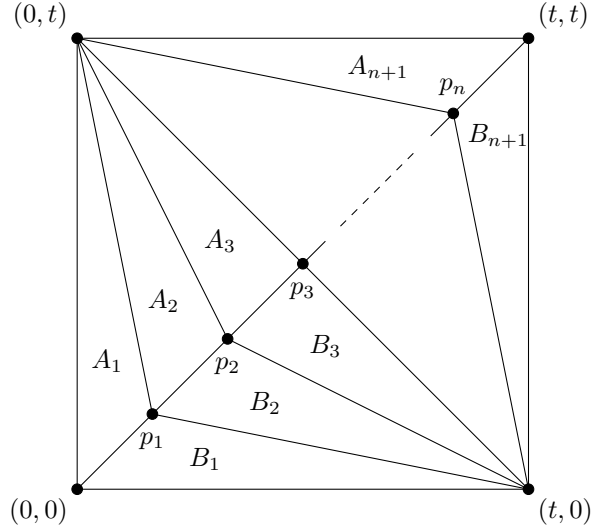
\begin{figure}
\begin{tikzpicture}
\draw (0,6) node[anchor=south east] {$(0,t)$} rectangle (6,0) node[anchor=north west] {$(t,0)$};
\draw[fill=black] (0,0) circle (2pt);
\draw[fill=black] (0,6) circle (2pt);
\draw[fill=black] (6,6) circle (2pt);
\draw[fill=black] (6,0) circle (2pt);
\draw (0,0) node[anchor=north east] {$(0,0)$}-- (3,3);
\draw (3,3) -- (3.3,3.3);
\draw[dashed] (3.5,3.5) -- (4.5,4.5);
\draw (4.7,4.7) -- (5,5);
\draw (5,5) -- (6,6) node[anchor=south west] {$(t,t)$};
\draw (0,6) -- (1,1) node[label=below:$p_1$] {};
\draw (0,6) -- (2,2) node[label=below:$p_2$] {};
\draw (0,6) -- (3,3) node[label=below:$p_3$] {};
\draw (0,6) -- (5,5) node[label=above:$p_n$] {};
\draw (6,0) -- (1,1);
\draw (6,0) -- (2,2);
\draw (6,0) -- (3,3);
\draw (6,0) -- (5,5);
\draw[fill=black] (1,1) circle (2pt);
\draw[fill=black] (2,2) circle (2pt);
\draw[fill=black] (3,3) circle (2pt);
\draw[fill=black] (5,5) circle (2pt);
\draw (1.7,0.4) node{$B_1$};
\draw (2.5,1.15) node{$B_2$};
\draw (3.3,1.9) node{$B_3$};
\draw (5.6,4.7) node{$B_{n+1}$};
\draw (0.4,1.7) node{$A_1$};
\draw (1.15,2.5) node{$A_2$};
\draw (1.9,3.3) node{$A_3$};
\draw (4,5.6) node{$A_{n+1}$};
\end{tikzpicture}
\caption{The triangulation $T_n$.}
\label{fig:diagonal}
\end{figure}
\end{center}

We will compute the polynomial $p(\T_n)$ for each $n$ using a straightforward
elimination of variables.  In particular, we will show that $\deg(\T_n)=n$.

For convenience we illustrate with drawings of a square in the plane with 
coordinates $(0,0)$, $(t,0)$, $(t,t)$, and $(0,t)$.  
We label the $n$ interior points $p_i=(x_i,y_i)$ (in ``increasing order''
as in Figure \ref{fig:diagonal}), and let $p_0=(0,0)$ and $p_{n+1}=(t,t)$.  Then $\T_n$
has edges connecting $p_i$ to
$p_{i+1}$, for $i=0,1,\ldots,n$, and also connecting each
$p_i$ to both corners $(0,t)$ and $(t,0)$.

For the areas we use the following notation:  for $1\leq i \leq n+1$, 
let $A_i$ denote twice the area of the 
triangle with vertices $\{(0,t), p_{i-1}, p_i\}$ and let $B_i$ denote 
twice the area of the triangle with vertices $\{(t,0), p_i, p_{i-1}\}$.
This accounts for all $2n+2$ triangles.

It is convenient to define the linear forms
\begin{eqnarray}\label{sigma}
\s_k := \sum_{i=1}^k (A_i+B_i), \\
\label{sigmabar}
\sb k := \sum_{i=k+1}^{n+1} (A_i+B_i) = 2t^2 - \s_k.
\end{eqnarray}
We allow $0\leq k\leq n+1$, i.e., $\s_0=0=\sb {n+1} $ and 
$\sb 0 =2t^2=\s_{n+1}$.

We now begin the elimination of variables.
In other words, beginning with the formulas expressing the
$A_i$ and $B_i$ in terms of $x_i,y_i,$ and $t$, we systematically
eliminate $x_i,y_i$, and $t$.

Here is the starting point:
$$
\begin{array}{rclcl}
A_i&= &\left|
\begin{array}{ccc} 1 & 1 & 1 \\
			0 & x_{i-1} & x_i \\
			t & y_{i-1} & y_i \\ 
			\end{array}
	\right|
&= &x_{i-1}y_i-x_iy_{i-1}+t(x_i-x_{i-1}); \\
\\
B_i&= &\left|
\begin{array}{ccc} 1 & 1 & 1 \\
			t & x_i & x_{i-1} \\
			0 & y_i & y_{i-1} \\ 
			\end{array}
	\right|
&= &-x_{i-1}y_i+x_iy_{i-1}+t(y_i-y_{i-1}).
\end{array}
$$

\begin{lem}
For each $0\leq k\leq n+1$, 
\begin{eqnarray}
y_k=\frac {\s_k} t  - x_k.
\end{eqnarray} 
\end{lem}

\begin{proof}  It is clear both from the above formulas and
from the picture (imagining a segment drawn from the origin
to $p_k$) that $t(x_k+y_k)=\s_k.$
\end{proof}

Using the lemma, we easily eliminate all $y$'s.  So, on to the $x$'s.
Substituting for the $y$'s in the formula for $B_k$:

\begin{eqnarray*}
B_k &= &-x_{k-1}y_k+x_ky_{k-1}+t(y_k-y_{k-1}) \\
&= & -x_{k-1}(\frac {\s_k} t - x_k) + x_k ( \frac {\s_{k-1}} t -x_{k-1})
+ (\s_k-\s_{k-1}) -tx_k+tx_{k-1} \\
&= & x_k(\frac {\s_{k-1}} t - t) -x_{k-1}(\frac {\s_k} t -t) + A_k + B_k,
\end{eqnarray*}
where on the last line we used the fact that $\s_k-\s_{k-1}=A_k+B_k$.
Solving for $x_k$ gives:
$$
x_k=\frac 1 t \frac { (tx_{k-1})(\s_k-t^2) - t^2 A_k }{(\s_{k-1}-t^2)}.
$$
(We have deliberately chosen to write a factor of $1/t$ in front,
even though it could be canceled.)  Now, by \eqref{sigma} 
we rewrite $\s_i-t^2=(\s_i -\sb i )/2$ and $t^2=\sb 0 /2=-(\s_0 -\sb 0 )/2$,
thus obtaining
\begin{eqnarray*}
x_k=\frac 1 t \frac {(tx_{k-1})(\s_k -\sb k ) + (\s_0 -\sb 0 )A_k } {(\s_{k-1} -\sb {k-1} )}.
\end{eqnarray*}

We now have $x_k$ written in terms of $x_{k-1}$, so we may 
recursively eliminate the $x$'s.

For shorthand we will use the notation
$\l_k $ for the linear form $\s_k - \sb k $,
and 
\begin{eqnarray}\label{lk}
\l_{k_1\cdots k_j} = \prod_{i=1}^j \l_{k_i}=\prod_{i=1}^j (\s_{k_i} - \sb {k_i} ).
\end{eqnarray}

So our expression for $x_k$ becomes
\begin{eqnarray}\label{xk}
x_k=\frac 1 t \frac {(tx_{k-1}) \l_k + A_k \l_0}{\l_{k-1}}.
\end{eqnarray}

\begin{lem}
For $1\leq k\leq n+1$, we have
$$x_k = \frac 1 t \frac {p_k}{q_k},$$
where $p_k$ is the degree $k$ polynomial
$$p_k=\l_{01\cdots k}\left(\frac {A_1}{\l_{01}} + \frac {A_2}{\l_{12}}
+\cdots + \frac {A_k}{\l_{k-1,k}}\right)
$$
and $q_k$ is the degree $k-1$ polynomial
$q_k=\l_{12\cdots(k-1)}$.  (In case $k=1$, $q_1$ is the empty product 1.)
\end{lem}

\begin{proof}
When $k=1$ the claim is that $x_1=A_1/t$, which is clearly true.  The rest follows 
by induction from 
(\ref{xk}).  Specifically, if $x_{k-1}=\frac 1 t \frac {p_{k-1}} {q_{k-1}}$
then plugging into (\ref{xk}) gives the recursions
$$q_k=q_{k-1}\l_{k-1}$$ and
$$p_k=p_{k-1}\l_k + A_k \l_0 q_{k-1}.$$
The formula for $q_k$ follows immediately, and the formula for $p_k$
is easily verified.
\end{proof}

The preceding lemma is valid for $k=n+1$, so we now have an explicit  
expression for $x_{n+1}$ in terms of $t$ and the $A_i$ and $B_i$.
But $x_{n+1}=t.$  This gives the relation
$t=x_{n+1}=\frac 1 t \frac {p_{n+1}}{q_{n+1}}$, which we write as
$-\l_0 q_{n+1} = 2p_{n+1}$ (recalling that $-\l_0=2t^2$).
Plugging in for $q_{n+1}$ and $p_{n+1}$ from the preceding
lemma gives
\begin{eqnarray}\label{eq:diagonal}
-\l_{01\cdots n} &= &2 \l_{01\cdots (n+1)}
\left(\frac {A_1}{\l_{01}} + \frac {A_2}{\l_{12}}
+\cdots + \frac {A_{n+1}}{\l_{n,n+1}}\right)\\
& = & 2\left(A_1 \l_{23\cdots (n+1)} + A_2 \l_{34\cdots n+1,0} +
\cdots + A_{n+1} \l_{01\cdots(n-1)} \right).
\end{eqnarray}

This is a homogeneous degree $n+1$ polynomial relation 
among the areas $A_i, B_i$, and the coefficients are integers with no
common factor.  Thus, by Theorem \ref{cor:existence} the polynomial
$p(\T_n)$ must be an irreducible factor of \eqref{eq:diagonal}.
Note that \eqref{eq:diagonal} itself is not 
irreducible:  both sides of the above equation contain a (linear) 
factor of $-\l_0=\l_{n+1}=\sum (A_i+B_i)$.  Since there exist nondegenerate 
drawings of $T_n$ (i.e., drawings in which $t\ne0$ and hence the total area
is nonzero), this linear factor is not equal to $p(\T_n)$.  If we factor this out, 
we are left with a degree $n$ polynomial; we next argue that this must be 
equal to $p(\T_n)$.  

One way to write the degree $n$ relation is to cancel $-\l_0$ from the left side against 
$\l_{n+1}$ from every term on the right except the last, and $\l_0$
from the remaining term.  This puts the remaining polynomial in the form
recorded in the following theorem.

\begin{thm}
Let $n\ge 1$ and let $\T_n=(T_n,\emptyset)$, where $T_n$ is the triangulation 
illustrated in Figure \ref{fig:diagonal}.  The polynomial $p(\T_n)$ is equal to 
\begin{equation}\label{poly}
\l_{12\cdots n}-2\left(
A_1 \l_{23\cdots n} + A_2 \l_{34\cdots n0} + \cdots + A_n\l_{01\cdots(n-2)}
- A_{n+1}\l_{12\cdots(n-1)}
\right).
\end{equation}
In particular, $\deg(\T_n)=n$.
\end{thm}

(See Equations \eqref{sigma}, \eqref{sigmabar}, and \eqref{lk} for the notation.)
Note that the polynomial \eqref{poly} has integer coefficients.  By modifying CoCoA
code originally written by D.~Perkinson, we computed $p(\T_n)$ for $n=1,2,3,4,5$,
and the CoCoA script produced this polynomial divided by $2^{n-1}$.
For $n\ge 6$ the problem was too large for our laptops.

\begin{proof}
We have shown that the areas $A_i$, $B_i$ satisfy \eqref{poly}.  If one were to prove
directly that this polynomial is irreducible, then it would follow that it equals $p(\T_n)$.
However, this is unnecessary, because Corollary \ref{cor:linearlowerbound} implies
that the degree of $p(\T_n)$ is at least $n$, and so \eqref{poly} must be the polynomial 
$p(\T_n)$ (and \eqref{poly} is irreducible).
\end{proof}

\end{example}

The next family of examples illustrates two notable phenomena.  First, we show
that $\deg(\T)$ can grow exponentially in the number of vertices of the triangulation.  
In a slightly different context, this was also observed in \cite{gonewild}.
In that paper, the authors give a general upper bound $\deg(\T) \le (t+1)2^t$
where $t$ is the number of non-degenerate triangles in $\T$.

Second, we will see that when running Algorithm \ref{alg}, it is possible to
get arbitrarily large multiplicities of components.

\begin{example}
We define triangulations $T_{n,k}$ for all non-negative integers $n$ and $k$ with 
$k\le \frac{n+1}{2}$.  Roughly speaking, we form $T_{n,k}$ by starting with the diagonal
case $T_n$ and ``exploding'' exactly $k$ length 2 segments along the diagonal into
pairs of triangles.  (In particular $T_{n,0}=T_n$.)  Precisely, $T_{n,k}$ has
a total of $n+k+4$ vertices:  the four corners $P=p_0,Q,R=p_{n+1},S$; 
the $n$ vertices $p_1,\ldots,p_n$ arranged diagonally as before; and also ``duplicates''
$p_1', p_3', \ldots, p_{2k-1}'$ of $k$ of the diagonal vertices.  The $2n+2k+2$ triangles of
$T_{n,k}$ have the following vertex sets (see Figure \ref{fig:exploded}):
\begin{itemize}
\item
$\{S,p_i,p_{i+1}\}$ for $i=0,\ldots,n$;
\item
$\{Q,p_{2i-1}',p_{2i-2}\}$ for $i=1,\ldots,k$;
\item
$\{Q,p_{2i},p_{2i-1}\}$ for $i=1,\ldots,k$;
\item
$\{Q,p_{i+1},p_i\}$ for $i=2k,\ldots,n$;
\item
$\{p_{2i-2},p_{2i-1}',p_{2i-1}\}$ for $i=1,\ldots,k$;
\item
$\{p_{2i},p_{2i-1},p_{2i-1}'\}$ for $i=1,\ldots,k$.
\end{itemize}

Note that one obtains $T_n$ from $T_{n,k}$ by contracting each
of the edges $p_ip_i'$.

\begin{center}
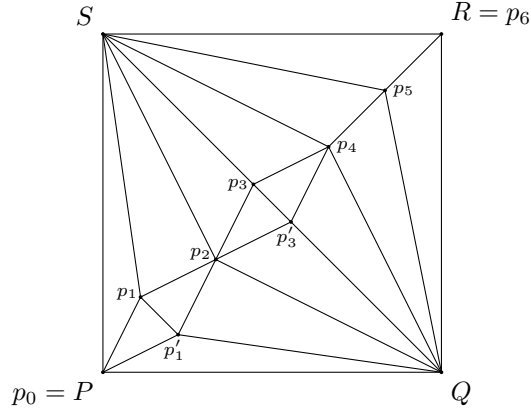
\begin{figure}
\begin{tikzpicture}[scale=.25]

\draw (0,18) node[anchor=south east] {$S$} 
	rectangle (18,0) node[anchor=north west] {$Q$}
	;
\draw[fill=black] (0,0) circle (2pt);
\draw[fill=black] (0,18) circle (2pt);
\draw[fill=black] (18,18) circle (2pt);
\draw[fill=black] (18,0) circle (2pt);
\draw (0,0) node[anchor=north east] {$p_0=P$} --(2,4) --(6,6) --(4,2) --(0,0) ;
\draw (2,4)--(4,2);
\draw (6,6)--(8,10)--(12,12)--(10,8)--(6,6);
\draw (8,10)--(10,8);
\draw (12,12)--(18,18) node[anchor=south west] {$R=p_6$};
\draw (18,0) -- (4,2);
\draw (18,0) -- (6,6);
\draw (18,0) -- (10,8);
\draw (18,0) -- (12,12);
\draw (18,0) -- (15,15);
\draw (0,18) -- (2,4);
\draw (0,18) -- (6,6);
\draw (0,18) -- (8,10);
\draw (0,18) -- (12,12);
\draw (0,18) -- (15,15);
\draw[fill=black] (2,4) circle (2pt);
\draw[fill=black] (4,2) circle (2pt);
\draw[fill=black] (6,6) circle (2pt);
\draw[fill=black] (8,10) circle (2pt);
\draw[fill=black] (10,8) circle (2pt);
\draw[fill=black] (12,12) circle (2pt);
\draw[fill=black] (15,15) circle (2pt);
\draw (1.3,4.2) node{{\scriptsize $p_1$}};
\draw (5.2,6.3) node{{\scriptsize $p_2$}};
\draw (3.8,1.2) node{{\scriptsize $p_1'$}};
\draw (7.2,10) node{{\scriptsize $p_3$}};
\draw (13,12) node{{\scriptsize $p_4$}};
\draw (9.8,7.2) node{{\scriptsize $p_3'$}};
\draw (16,15) node{{\scriptsize $p_5$}};
\end{tikzpicture}
\caption{The exploded diagonal triangulation $T_{5,2}$.}
\label{fig:exploded}

\end{figure}

\end{center}

\begin{thm}\label{thm:exploded}
For the triangulation $\T_{n,k}=(T_{n,k},\emptyset)$, we have 
$$\deg(\T_{n,k}) \ge 2^{k-1}(2n-k).$$
\end{thm}

\begin{proof}
We proceed by induction on $k$, by running the algorithm.  If
$k=0$, the triangulation is the diagonal case and the result agrees
with the previous example.  

Suppose $k>0$.
We first kill the triangle $\{p_{2k},p_{2k-1},p_{2k-1}'\}$.  The result
has a single component.  We then kill the triangle $\{p_{2k-2},p_{2k-1}',p_{2k-1}\}$.
The result has two components.  One is obtained by contracting the edge
$p_{2k-1}p_{2k-1}'$; this triangulation is exactly $T_{n,k-1}$.

The other component, $\dot\T$, has a condition consisting of the four vertices
$\{p_{2k-2},p_{2k-1},p_{2k-1}',p_{2k}\}$.  We analyze this using the following 
\begin{ut}
Suppose $\T=(T,C)$ is a triangulation, where $C\supset\{w,x,y,z\}$ and
$wxy$ and $xyz$ are triangles of $T$.
Let $\hat T$ be the triangulation
obtained from $T$ by replacing the triangles $wxy$ and $xyz$ with triangles
$wxz$ and $wyz$, and let $\hat\T=(\hat T,C)$.  Then $p(\T)=p(\hat \T)$.
\end{ut}

The trick works because $X(\T)=X(\hat \T)$, $Y(\T)=Y(\hat \T)$, and $f_\T=f_{\hat \T}$.
(Note that $f_{(T,\emptyset)} \ne f_{(\hat T,\emptyset)}$.)  Thus $V_\T=V_{\hat\T}$
and $p_\T=p_{\hat\T}$.

In this case, in $\dot T$ we exchange the edge $p_{2k-1}p_{2k-1}'$ for the
edge $p_{2k-2}p_{2k}$, and the resulting triangulation is a subdivision of 
$T_{n-1,k-1}$.  Thus $\deg(\dot\T)=\deg(\T_{n-1,k-1})$ and so we have
$$\deg(\T_{n,k}) \ge \deg(\T_{n,k-1}) + \deg(\T_{n-1,k-1}).$$
Since we know $\deg(\T_{n,0})=n$, the theorem follows by solving the above
recurrence.
\end{proof}

Theorem \ref{thm:exploded} implies the following table of lower bounds, 
with answers known to be correct in bold.

\medskip

\begin{center}
\begin{table}[h]
\begin{tabular}{r|cccccccc}
$n$ & $0$ & $1$ & $2$ &$3$&$4$&$5$&$6$&$7$\\
\hline \\
$k=0$ & ${\bf 1}$ & ${\bf 1}$ & ${\bf 2}$ & ${\bf 3}$ & ${\bf 4}$ & ${\bf 5}$ & ${\bf 6}$ & ${\bf 7}$ \\
$1$ & & ${\bf 2}$ & ${\bf 3}$ & $5$ & $7$ & $9$ & $11$ & $13$ \\
 $2$ & & & & $8$ & $12$ & $16$ & $20$ & $24$ \\
 $3$ & & & & & & $28$ & $36$ & $44$ \\
 $4$ & \phantom{80}& \phantom{80} &  \phantom{80}&  \phantom{80}& \phantom{80} & \phantom{80} & \phantom{80} & $80$ \\
 \\
\end{tabular}

\caption{Lower bounds for $\deg(\T_{n,k})$.  Values in bold are known to be sharp.}
\label{table}
\end{table}
\end{center}

Consider for a moment triangulations $\T=(T,C)$ with no subdivisions, $C=\emptyset$,
and $k$ interior vertices.
Example \ref{example1} is the only such triangulation with $k=1$,
and of course its polynomial is linear.
Up to isomorphism, Example \ref{example2} is the only such triangulation with $k=2$,
and it has degree 2.  There are exactly two isomorphism classes of these triangulations 
with $k=3$, namely $\T_3$ and $\T_{2,1}$; both have degree 3.  The triangulation
$\T_{3,1}$ is therefore the smallest triangulation without subdivisions satisfying 
$\deg(\T)>k$.

We close with an observation about the multiplicities of components arising
in Algorithm \ref{alg}.  The above table contains the correct degrees if each 
component has multiplicity 1, when we run the algorithm as described.
However, there are other ways to run the algorithm, and for example if
$k \le \frac{n-1}2$ then we could instead begin with $\T_{n,k}$ and kill triangles
$Qp_np_{n-1}$ and then $Sp_{n-1}p_n$.  This would produce the recurrence
$$\deg(\T_{n,k}) \ge \deg(\T_{n-1,k}) +1.$$
Now, if Table \ref{table} contains the correct values for $\deg(\T_{n,k})$
for all $n$ and $k$, then for fixed $k$ and $n\ge 2k+1$, the linear factor arising
from this factorization must occur with multiplicity $2^k$.  In view of the discussion
after Conjecture \ref{conj:algworks}, this means that the geometric intersection of
$V$ with the corresponding hyperplanes consists entirely of singular points of $V$.


\end{example}

\bibliography{../..//biblio.bib}
\bibliographystyle{plain}

\end{document}